\documentclass[12pt]{article}

\usepackage{amsmath, amsthm, amssymb, amsfonts, enumerate}
\usepackage[applemac]{inputenc} 
 \usepackage[french] {babel}

\def\dis{\displaystyle}
\def \R{I\!\!R}
\def \E{I\!\!E}
\def \P{I\!\!P}
\def \N{I\!\!N}

\newtheorem{thm}{Theorem}[section]
\newtheorem{cor}[thm]{Corollary}

\newtheorem{lem}[thm]{Lemma}
\newtheorem{pro}[thm]{Proposition}
\newtheorem{defi}[thm]{Definition}
\newtheorem{rem}[thm]{Remark}
\newtheorem{nota}[thm]{Notation}
\newtheorem{notas}[thm]{Notations}

\newtheorem{cla}[thm]{Claim}

\def \lim{ {\rm lim} }

\def \inf{ {\rm inf} }
\def \sup{ {\rm sup} }

\def\abstract{\begin{center} \small\bf Abstract\end{center}\small}

\title{The value of   Repeated Games  
with an informed controller}
 
\author{J\'er\^ome Renault\thanks {CMAP and Economic Department, Ecole Polytechnique, 91128 Palaiseau Cedex, France.  email: jerome.renault@polytechnique.edu}}

\date{\today} 

\begin{document}     

\maketitle
\begin{abstract}  We consider the  general model of zero-sum repeated games (or stochastic games with signals), and assume that one of the players is fully  informed and controls the transitions of the state variable. We prove the existence of the uniform value, generalizing several results of the literature. A preliminary existence result is obtained for a certain class of  stochastic games played with pure strategies.   \\


\end{abstract}


\section{Introduction}

The context of this work is the characterization of      repeated game models where the  value exists. We first consider here  general    repeated games defined with finite sets of states, actions and signals.  They contain   usual  stochastic games,  standard repeated games with incomplete information and also repeated games with signals. At each stage the players will play a matrix game depending on a  parameter called state, this   state is partially known and evolves  from stage to stage,  and after each stage every player receives some  private signal on the current situation.  We make two   important hypotheses. We first  assume that player 1 is informed, in the sense that he can always deduce the current state and player 2's signal from his own signal. Secondly, we assume that player 1 controls the transitions, in the sense that the law of the couple (new state, signal received by player 2) does not depend on player 2's actions. We call  ``repeated games with an informed controller" the  games satisfying these two hypotheses.

This class of games  contains  Markov chain repeated games with lack of information on one side as studied in Renault, 2006, and is  more general since, in particular,  we allow here for transitions of the state depending also on player 1's actions. It  also contains  stochastic games with a single controller and incomplete information on the side of his opponent, as studied in Rosenberg {\it et al.}, 2004 (see subsection \ref{subsecappli} here).  And {\it a fortiori} it contains the   standard     repeated games with incomplete information on one side and perfect monitoring introduced by  Aumann and Maschler. Notice that   repeated games with an informed controller   contains (weak) forms of the three  main aspects of   general repeated games models: the stochastic aspect (the state evolves from one stage to another and is controlled here by player 1), the incomplete information aspect (player 2   has an incomplete knowledge of the state), and the signalling aspect (players observe signals rather than actions).  And we believe that  the existence result presented here is the first one to significantly deal   with   these three  aspects simultaneously. On the contrary, they do not contain  stochastic games, where the transitions are controlled by both players (see Mertens Neyman 1981 for the existence of the uniform value in such games).

We   prove  the existence of the uniform value via  several steps, and several games are considered. These  are:  our original repeated game where player 1 is informed and controls the transitions (level 1), an auxiliary stochastic game (level 2), and finally a one-player repeated game, i.e. a dynamic programming problem (level 3).  A  crucial point is that in our original game, the set of states  $K$ is finite. 

The  auxiliary stochastic game has  the following features. It is played with pure strategies. The new set of states is $X=\Delta(K)$, the set of probabilities over $K$, and represents in the original game the belief of player 2 on the current state\footnote{The idea of considering   an auxiliary stochastic game is certainly  not new, see for example Mertens 1986, Coulomb, 2003 or Mertens {\it et al.}, 1994, Part A, Ch IV, section 3.}. In the auxiliary stochastic game, the new state is known to both players, and actions played are perfectly observed after each stage. It is also convenient  to consider    states in  $\Delta_f(X)$, the set of probabilities with finite support on $X$. To express an informational gap, we use  the Choquet order of sweeping of probabilities:   given $u$ and $v$  in $\Delta_f(X)$, we say that  $u$ is better than $v$, or $v$  is a sweeping of $u$, if for every concave continuous mapping $f$ from $X$ to the reals, $u(f)\geq v(f)$. And it is essentially possible  to model the original informational advantage of player 1  with a ``splitting hypothesis" defined via this order. For the topological part, we use the weak* topology on the set $\Delta(X)$ of Borel probabilities on $X$, and more precisely the Wasserstein distance. This allows to carefully control the Lipschitz constant of the value functions.

A dynamic programming problem can be derived from the auxiliary stochastic game. The role of player 2 disappears, and the set of states of the dynamic programming problem is $Z=\Delta_f(X)\times [0,1]$. $\Delta_f(X)$ is dense in $\Delta(X)$ for the weak* topology, so $Z$ can be viewed as a precompact metric space.  We define, for every $m$ and $n$, a value $w_{m,n}$ as the supremum payoff player 1 can achieve when his payoff is  defined as the {\it minimum}, for $t$ in $\{1,..,n\}$, of his  average rewards computed between stages $m+1$ and $m+t$.  It is possible to prove that the family $(w_{m,n})$ is uniformly equicontinuous, and together with the precompacity of the state space this  implies the existence of the uniform value for the  the dynamic programming problem. The proof of this implication can be found in a  companion paper (see Corollary 3.8, Renault 2007), which only deals with 1-player games and  can be read independently.

The present paper is organized as follows. In section \ref{secstoch}, we consider a particular class  of stochastic games including our auxiliary games of level 2.  We think this class of games may  be interesting in itself. It is defined with hypotheses making no reference to the original finite set $K$, and  we prefer to start by presenting  this class, which can  be considered as  both more general and simpler to study than the auxiliary stochastic games. We  prove that these games have a uniform value using the  result on dynamic programming proved in Renault, 2007. In section \ref{sec45}, we consider our original repeated game and show how the existence of the uniform value in these games is implied by the existence of the uniform value for the stochastic games of section \ref{secstoch}. Finally, we obtain    formulas expressing  the uniform value in terms of the values of some finite games. More precisely, let  $v_{m,n}$ be  the value of the game where the global  payoff is  defined as the   average of the   payoffs between stage $m+1$ and stage $m+n$.   We show  in particular that $\inf_{n\geq 1} \sup_{m\geq 0}  v_{m,n}= \sup_{m\geq 0}\inf_{n\geq 1} v_{m,n}$, and this is the uniform value (see subsection \ref{subsub11}). We conclude by discussing several hypotheses and present a few open problems.

\section{A certain class of  stochastic games} \label{secstoch}

\subsection{Model}\label{sub29}
 
We consider in this section 2-player zero-sum stochastic games with complete information and standard observation,  {\it played with pure strategies}. We assume that after each stage, the new state is selected according to a probability with finite support. 

 If $X$ is a non empty set, we denote by $\Delta_f(X)$ the set of probabilities on $X$ with finite support. 
 
 We consider: 

 \begin{quote}

 $\bullet$ three non empty  sets:  a set of states $X$,    a set $A$ of actions for player 1, and  a set $B$ of actions   for player 2,  
 
 $\bullet$ an element $u$ in  $\Delta_f(X)$, called the initial  distribution on states,  
  
$\bullet$  a mapping $g$ from $X \times A \times B$ to $[0,1]$, called the payoff function of player 1, and
 
 $\bullet$ a mapping $l$ from $X\times A \times B$ to   $\Delta_f(X)$, called the transition function.  

 \end{quote}
 
\noindent    The interpretation  is the following. The initial state $p_1$ in $X$ is selected according to $u$, and is announced to both players. Then simultaneously, player 1 chooses $a_1$ in $A$, and player 2 chooses $b_1$ in $B$. The stage payoff is $g(p_1,a_1,b_1)$ for player 1, and $-g(p_1,a_1,b_1)$ for player 2, then  $a_1$ and $b_1$ are publicly observed, and a new state $p_2$ is selected according to   $l(p_1,a_1,b_1)$, etc... At any stage $t\geq2$, the state $p_t$ is selected according to $l(p_{t-1},a_{t-1},b_{t-1})$, and announced to both players. Simultaneously, player 1 chooses $a_t$ in $A$ and player 2 chooses $b_t$ in $B$. The stage payoffs are $g(p_t,a_t,b_t)$ for player 1 and the opposite for player 2. Then $a_t$ and $b_t$ are publicly announced, and the play proceeds to stage $t+1$. \\

From now on we fix $\Gamma=(X,A,B,g,l)$, and for every $u$ in $\Delta_f(X)$ we denote by $\Gamma(u)=(X,A,B,g,l, u)$ the corresponding stochastic game induced by $u$. For the moment we make no assumption on $\Gamma$. We   start with elementary definitions and notations. 

A strategy for player 1 is a sequence $\sigma=(\sigma_n)_{n \geq 1}$, where for each $n$, $\sigma_n$ is a mapping from $(X \times A \times B)^{n-1}\times X$ to $A$, with the interpretation that $\sigma_n(p_1,a_1,b_1,...,p_{n-1},a_{n-1},b_{n-1},p_n)$ is the action played by player 1 at stage $n$ after $(p_1,a_1,b_1,...,p_{n-1},a_{n-1},b_{n-1},p_n)$ occurred. $\sigma_1$ is just a mapping from $X$ to $A$ giving the first action played by player 1 depending on the initial state. Similarly, a strategy for player 2 is a  sequence $\tau=(\tau_n)_{n \geq 1}$, where for each $n$, $\tau_n$ is a mapping from $(X \times A \times B)^{n-1}\times X$ to $B$. We denote by $\Sigma$ and ${\cal T}$ the sets of strategies of player 1 and player 2, respectively.

Fix for a while $(u,\sigma, \tau)$, and assume that player 1 plays $\sigma$ whereas player 2 plays $\tau$ in the game $\Gamma(u)$. The initial state $p_1$ is selected according to $u$, then the first actions are $a_1=\sigma_1(p_1)$ and $b_1=\tau_1(p_1)$. $p_2$ is selected according to $l(p_1,a_1,b_1)$, then $a_2=\sigma_2(p_1,a_1,b_1,p_2)$, $b_2=\tau_2(p_1,a_1,b_1,p_2)$, etc... By induction this defines, for every positive $N$, a probability with finite support on the set $(X\times A \times B)^N$ corresponding to the set of the first $N$ states and actions.  It is standard   that these probabilities can be uniquely extended to a probability $\P_{u,\sigma, \tau} $ on the set of plays $\Omega=(X\times A \times B)^{\infty}$, endowed with the $\sigma$-algebra generated by the cylinders (one can apply, e.g., theorem 2.7.2. p.109 in Ash, 1972).  

\begin{defi}\label{def5}   The average expected payoff for player 1
 induced by $(\sigma,
\tau)$ at the first $N$ stages in the game $\Gamma(u)$ is denoted by: 
$$\gamma_N^u (\sigma, \tau)=\E_{\P_{u, \sigma, \tau}} \left( {1 \over N}\sum_{n=1}^N g (p_n,a_n,b_n)
\right).$$  \end{defi}

\begin{defi}\label{def6}   For $u$ in $\Delta_f(X)$ and $N\geq 1$,   the game  $\Gamma_N(u)$ is the zero-sum game with normal form $(\Sigma, {\cal T}, \gamma_N^{u})$. \end{defi}

\noindent  $\Gamma_N(u)$ is called the $N$-stage  game with initial distribution $u$,  and corresponds to  the one-shot game where player 1's  strategy set is $\Sigma$, player 2's strategy set is ${\cal T}$, and $\gamma_N^u$ is the payoff function for player 1. It has a value if:
$\sup_{\sigma \in \Sigma} \inf_{\tau \in {\cal T}}\gamma_N^{u}(\sigma, \tau)= \inf_{\tau \in {\cal T}}\sup_{\sigma \in \Sigma}\gamma_N^{u}(\sigma, \tau).$ A strategy $\sigma$, if any,  achieving the supremum on the LHS   is then called an optimal strategy for player 1. Similarly,  a strategy $\tau$, if any,  achieving the infimum  on the RHS is then called an optimal strategy for player 2.

\begin{notas}\label{nota1} $\;$

 \noindent For $p$ in $X$, we denote by $\delta_p \in\Delta_f(X)$ the Dirac measure on $p$.  A probability $u$ in $\Delta_f(X)$ is written $u=\sum_{p \in X} u(p) \delta_p$, where $u(p)$ is the probability of $p$ under $u$.    
 \vspace{0,5cm}

 \noindent  For $u$ in $\Delta_f(X)$, if  $\Gamma_N(u)$ has a value, we denote it by $\tilde{v}_N(u)\in [0,1]$. For  $p$ in $X$,  if   $\Gamma_N(\delta_p)$ has a value, we denote it by $v_N(p)$ and we have $v_N(p)=\tilde{v}_N(\delta_p)$.\end{notas}

\noindent Notice that if $p\neq p'$, then $\delta_{1/2p+1/2p'}\neq 1/2 \, \delta_{p} +1/2 \, \delta_{p'}$, so we will not identify a state $p$ with the measure $\delta_p$. When the value of the $N$-stage game exists for every initial distribution, $\tilde{v}_N$ is a mapping from $\Delta_f(X)$ to $\R$, whereas $v_N$ is a mapping from $X$ to $\R$. It is easy to see that:  $\forall u \in \Delta_f(X), \forall N \geq 1, \forall \sigma \in \Sigma, \forall \tau \in {\cal T},$
$$\P_{u, \sigma, \tau}= \sum_{p \in X} u(p) \P_{ \delta_p, \sigma, \tau} \;\; {\rm and }\;\; \gamma_N^u (\sigma, \tau)=\sum_{p \in X} u(p) \gamma_N^{\delta_p}(\sigma, \tau).$$
\begin{cla} \label{cla2,8} If $v_N(p)$ exists for each $p$ in $X$, then $\tilde{v}_N(u)$ exists for every $u$ in $\Delta_f(X)$ and $\tilde{v}_N(u)=\sum_{p \in X}u(p)v_N(p)$.\end{cla}


\vspace{0,3cm}

We now consider an  infinite time horizon.
 
 \begin{defi}\label{def7} Let $u$ be in $\Delta_f(X)$.

The lower (or maxmin) value of $\Gamma(u)$ is: $$\underline{v}(u) = \sup_{\sigma\in \Sigma} \liminf_n \left( \inf_{\tau \in {\cal T} }\gamma_n^u( \sigma, \tau)\right).$$

 The upper (or minmax) value of $\Gamma(u)$ is:
 
  $$\overline{v}(u) = \inf_{\tau \in {\cal T} }\limsup_n \left(\sup_{\sigma\in \Sigma}  \gamma_n^u( \sigma, \tau)\right).$$
 
\noindent $ \underline{v}(u)\leq \overline{v}(u)$.  $\Gamma(u)$ is said to have a uniform value if and only if $\underline{v}(u)=\overline{v}(u)$, and in this case the uniform value is $\underline{v}(u)=\overline{v}(u)$. \end{defi}

An equivalent definition of the uniform value is as follows. Given a real number $v$,  we say that player 1 can guarantee $v$ in $\Gamma(u)$  if: $\forall \varepsilon>0, \exists \sigma \in \Sigma, \exists N_0, \forall N\geq N_0, \forall \tau\in {\cal T}, \gamma_N^u(\sigma, \tau) \geq v- \varepsilon.$ Player 2 can guarantee $v$ in $\Gamma(u)$ if: $\forall \varepsilon>0, \exists \tau\in {\cal T}, \exists N_0, \forall N\geq N_0, \forall \sigma \in \Sigma, \gamma_N^u(\sigma, \tau) \leq v+ \varepsilon.$ If player 1 can guarantee $v$ and player 2 can  guarantee $w$ then clearly $w\geq v$. We also have: 
\begin{cla} \label{cla3} 

 $\underline{v}(u)=\max\{v \in \R, {\rm player \; 1\; can \; guarantee\; } v\; {\rm in} \;\Gamma(u)\;  \},$ 

\hspace{1,7cm}$\overline{v}(u)=\min\{v \in \R, {\rm player \;2 \; can \; guarantee \;} v\; {\rm in}\; \Gamma(u)\;\}.$

A real number  $v$ can be guaranteed by both players if and only if  $v$ is   the uniform value of $\Gamma(u)$.

\end{cla}

  Assume now that $\tilde{v}_N(u)$ exists for each $N$. If player 1 (resp. player 2)  can guarantee $v$ then $ \liminf_N\tilde{v}_N(u)\geq v$ (resp.  $\limsup_N\tilde{v}_N(u)\leq v$). As a consequence we have:
  
  \begin{cla}\label{cla4} Assume that $\tilde{v}_N(u)$ exists for each $N$.  $$\underline{v}(u) \leq \liminf_N\tilde{v}_N(u)\leq  \limsup_N\tilde{v}_N(u)\leq \overline{v}(u).$$
\end{cla}

\noindent So in this case, the existence of the uniform value $\underline{v}(u)=\overline{v}(u)$ implies the existence of the ``limit value" $\lim_N \tilde{v}_N(u)$, and all  notions coincide: $\underline{v}(u)=\overline{v}(u)=\lim_N \tilde{v}_N(u).$ 

\subsection{An existence result for the uniform value}

We are interested in the existence of the uniform value, and are now going to make some hypotheses on  $\Gamma$.

\begin{rem} \label{rem3,5} We have in mind, in view of   application in section \ref{sec45},  the case of a repeated game with lack of information on one side where player 1 is informed and controls the transition. In these games, there is an underlying finite set of parameters $K$, and $X$ is the set of probabilities over $K$. Initially, a parameter $k$ is selected according to $p$, and is announced to player 1 only. Then   the parameter may change from stage  to stage, but it is always known by player 1 and its evolution is independent of player 2's actions.  It will be possible to check the following hypotheses H1 to H7   in this model. Our point here is to  be  more general and simpler. We want to be able to write a model without any reference to the underlying set of parameters, and  where players use pure strategies. \end{rem}

We first make the very important assumption that player 1 only controls the transitions:\\

\noindent {\bf Hypothesis H1}: the transition  $l$ does not depend on player 2's actions, i.e. $\forall p \in X, \forall a \in A, \forall b \in B, \forall b'\in B, l(p,a,b)=l(p,a,b').$\\

In the sequel we consider $l$ as a mapping from $X \times A$ to $\Delta_f(X)$, and we write $l(p,a)$ for the distribution on the next state if the actual state is $p$ and player 1 plays  $a$.\\

\noindent {\bf Hypothesis H2}: $X$ is a compact convex subset of a  normed vector space. \\

We denote by  $\Delta(X)$ the set of Borel probability measures on $X$, and we  will use on $\Delta(X)$ the weak* topology and the   Choquet order. $\Delta_f(X)$ is now seen as a subset of $\Delta(X)$. We first  fix notations and recall some definitions. We start with the topological aspect :  $X$ is in particular a compact metric space, and we denote by $d(p,q)$ the  distance between two elements $p$ and $q$ of $X$.

\begin{notas}\label{nota1,5} $\;$  We denote by $E$ the set of continuous mappings from $X$ to the reals, and by $E_1$ the set of   non expansive (i.e. Lipschitz with constant 1) elements of $E$.   For $u$ in $\Delta(X)$ and $f$ in $E$ we write $\dis u(f)=\int_{p \in X}f(p)du(p)$. Given $f$ in $E$, we extend $f$ by duality  to an affine  mapping $\tilde{f}: \Delta(X)\longrightarrow \R$ by $\tilde{f}(u)=u(f)$.  \end{notas}

In the following, $\Delta(X)$ will {\it  always} be endowed with the weak* topology: a sequence $(u_n)_n$ converges to $u$ in $\Delta(X)$ if and only if $u_n(f)\longrightarrow_{n \to \infty} u(f)$ for every $f$ in $E$. $\Delta(X)$ is itself compact,  the weak* topology can be metrized, and the set $\Delta_f(X)$ of probabilities on $X$ with finite support is dense in $\Delta(X)$ (see for example  Doob, 1994, Ch.VIII, section 5, and Malliavin, 1995, p.99).

 \begin{rem} \label{rem4}  An important  distance on $\Delta(X)$ which metrizes the weak* topology is the  following  (Fortet-Mourier-)Wasserstein distance, defined by: $$\forall u \in \Delta(X), \forall v \in \Delta(X),  \;\; d(u,v)= \sup_{f\in E_1} |u(f)-v(f)|.$$ 
 One can check that this distance has the following nice properties.   For every $p$, $q$ in $X$, $d(p,q)=d(\delta_p, \delta_q)$. Moreover,  for  $f$ in $E$ and $C\geq 0$,   {$f$ is $C$-Lipschitz if and only if $\tilde{f}$ is $C$-Lipschitz.}\end{rem}
 
 \vspace{0,5cm}

We will also use the convexity of $X$. In zero-sum games with lack of information on the side of player 2, it is well known that the value is a concave function of the parameter $p$: this fundamental property represents the advantage for player 1 to be informed (see for  example, Sorin 2002, proposition 2.2 p. 16). In our setup, we want the initial distribution $\delta_{1/2p+1/2p'}$ to be more advantageous for player 1 than the initial distribution $1/2\, \delta_{p} +1/2 \, \delta_{p'}$.  This is perfectly represented by the following   order:

\begin{defi} \label{def8} For $u$ and $v$ in $\Delta(X)$, we say that $u$ is better than $v$, or that $v$ is a sweeping of $u$, and we write $u\succeq v$, if:  
 
\centerline{for every {\it concave}  mapping $f$ in $E$, $u(f)\geq v(f).$} \end{defi}

This order was introduced by Choquet\footnote{For convenience, we reverse here Choquet's order, i.e. we write $u\succeq v$ instead of $u\preceq v$. For $u$ and $v$ in $\Delta_f(X)$, a simple characterization of $u\succeq v$ will be stated later, see proposition \ref{pro4}.}(1960). It is actually an order on $\Delta(X)$,   the maximal elements are the Dirac measures, and Choquet   proved that the minimal elements are the measures with support in the set of extreme points of $X$ (see    P.A. Meyer, 1966, theorem 24 p. 282). For every $f$ in $E$, we easily  have the equivalence: 
\begin{cla} \label{cla6}{$f$ is concave if and only if $\tilde{f}$ is non decreasing.}\end{cla}

 We   now  define   hypotheses H3 to H7.\\

\noindent {\bf Hypothesis H3}: $A$ and $B$ are compact convex subsets of topological vector spaces. \\

\noindent {\bf Hypothesis H4}:  For every  $(p,b)$ in $X \times B$, $(a\longrightarrow g(p,a,b))$ is concave and upper semi-continuous.  For every $(p,a)$ in $X \times A$, $(b\longrightarrow g(p,a,b))$ is convex and lower semi-continuous.\\

We will    prove in the sequel a natural dynamic programming principle {\it a la } Shapley (or Bellmann). 
 
 \begin{nota} \label{nota2} For $f$ in $E$ and $\alpha$ in $[0,1]$, we define $\Phi(\alpha,f):X \longrightarrow \R$ with:
 $$\forall p \in X, \;\; \Phi(\alpha,f)(p)= \sup_{a \in A} \inf_{b \in B} \;\left( \;  \alpha \; g(p,a,b) +(1-\alpha) \; \tilde{f} (l(p,a))\; \right).$$ \end{nota}

 \noindent {\bf Hypothesis H5}: There exists a subset ${\cal D}$ of $E_1$ containing $\Phi(1,0)$, and such that $ \Phi(\alpha,f) \in {\cal D}$ for every $f$ in ${\cal D}$ and $\alpha$ in $[0,1]$.\\

  \noindent {\bf Hypothesis H6}: For every  $(p,b)$ in $X \times B$, $(a\mapsto l(p,a,b))$ is continuous and concave.\\
  
    \noindent {\bf Hypothesis H7}: ``Splitting" Consider a   convex combination $p=\sum_{s=1}^S \lambda_s p_s$ in the set of states $X$,  and   a family of actions $(a_s)_{s \in S}$ in $A^S$. Then there exists $a$ in $A$ such that:
     $$\dis  l(p,a) \succeq \sum_{s\in S} \lambda_s l(p_s,a_s) \; \;\; \rm { and }\it  \;\; \; \inf_{b \in B} g(p,a,b) \geq \sum_{s\in S} \lambda_s \inf_{b \in B} g(p_s,a_s,b).$$
 
H3 and H4 are standard and, by Sion's theorem, will lead to the existence of the value of the stage game.   H5 is very important and will ensure that all value functions are 1-Lipschitz. We will provide later  a simple condition  implying H5, see remark \ref{rem5,3} . H6 is the only hypothesis where the topology on $\Delta(X)$ appears, and does not depend on a particular distance metrizing the weak* topology. H7 is the generalization of the well known splitting lemma for games with lack of information on one side. Under the hypotheses H1,..., H7, our main result   in theorem \ref{thm3} will be the existence of the uniform value.  We will also obtain several other properties, which will be expressed via the  following notions. \\

\begin{defi} \label{def10} A strategy $\sigma=(\sigma_t)_{t \geq 1}$ of player 1 is Markov if for each stage $t$, $\sigma_t$ only depends on the current state $p_t$.  A Markov strategy for player 1 will  be seen as a sequence 
$\sigma=(\sigma_t)_{t\geq 1}$, where for each $t$ $\sigma_t$ is a mapping from $X$ to $A$ giving the action to be played on stage $t$ depending on the current state. We denote the set of Markov strategies for player 1 by: $\Sigma^{M}=\{\sigma=(\sigma_t)_{t \geq1}, \rm{ with } \; \forall t, \sigma_t:X \longrightarrow A\}.$ Markov strategies for player 2 are defined similarly.
\end{defi}

\begin{defi} \label{def10,5} For $m\geq 0$ and $n \geq 1$,  the average expected payoff for player 1
 induced by a strategy pair $(\sigma,\tau)$ in $\Sigma \times {\cal T}$ from    stage   $m+1$ to stage $m+n$ is  denoted by: 
$$\gamma_{m,n} ^{u} (\sigma, \tau) =\E_{\P_{u, \sigma, \tau}} \left(  \frac{1}{n}\sum_{t=m+1}^{m+n}  \; g(p_t,a_t,b_t) \right).$$  \end{defi}
 
\begin{thm} \label{thm3} Assume that H1,..., H7 hold. 

Then  for every initial distribution $u$, the game $\Gamma(u)$ has  a uniform value $v^*(u)$. \\

Every player can guarantee $v^*(u)$ with a Markov strategy:   $\forall \varepsilon>0, \exists \sigma \in \Sigma^M,  \exists \tau\in {\cal T}^M, \exists N_0, \forall N\geq N_0, \forall \sigma' \in \Sigma, \forall \tau'\in {\cal T},$   $\gamma_N^u(\sigma, \tau') \geq v- \varepsilon \;  
and \; \gamma_N^u(\sigma', \tau) \leq v+ \varepsilon.$
\begin{eqnarray*}
We\;  have:\;  v^*(u) & = &\inf_{n\geq1} \sup_{m\geq 0} \tilde{v}_{m,n}(u) = \sup_{m\geq 0}  \inf_{n\geq1} \tilde{v}_{m,n}(u)\\
& =& \inf_{n\geq1} \sup_{m\geq 0} w_{m,n}(u) = \sup_{m\geq 0}  \inf_{n\geq1} w_{m,n}(u),
\end{eqnarray*}
where $\tilde{v}_{m,n}(u)= \sup_{\sigma \in \Sigma} \inf_{\tau \in {\cal T}} \gamma_{m,n} ^{u} (\sigma, \tau)=\inf_{\tau \in {\cal T}}  \sup_{\sigma \in \Sigma}  \gamma_{m,n} ^{u} (\sigma, \tau)$, 
 and 
$w_{m,n}(u)= \sup_{\sigma \in \Sigma} \inf_{\tau \in {\cal T}} \min_{t\in\{1,...,n\}} \gamma_{m,t} ^{u} (\sigma, \tau).$

For every $m$ and $n$, $\tilde{v}_{m,n}$ and $w_{m,n}$ are non expansive, and $(v_n)_n$ uniformly converges to $v^*$.
\end{thm}

\subsection{Proof of theorem \ref{thm3}}

We now prove theorem \ref{thm3}, and assume that H1,..., H7 hold. In the proof, we endow $\Delta(X)$ with the Wasserstein distance.  We denote by $\N^*$ the set of positive integers. By H3 and H4, for every $(p,a)\in X \times A$, the infimum is achieved in $\inf_{b \in B} g(p,a,b)$, and we   will simply write $g(p,a)$ for $\min_{b \in B} g(p,a,b)$. For each $p$, $(a\longrightarrow g(p,a))$ still is concave and upper semi-continuous.

\begin{lem} \label{lem3,5} For every concave $f$ in $E$ and $\alpha$ in $[0,1]$, $\Phi(\alpha,f)$ is concave. \end{lem}

\noindent{\bf Proof:}   Fix a   convex combination $p=\sum_{s=1}^S \lambda_s p_s$ in   $X$, and consider for each $s$ an element  $a_s$ in $A$. By the splitting hypothesis H7, one can find $a$ in $A$ such that $  l(p,a) \succeq \sum_{s\in S} \lambda_s l(p_s,a_s) \; \;\; \rm { and }\it  \;\; \; g(p,a) \geq \sum_{s\in S} \lambda_s   g(p_s,a_s).$ $f$ is concave so $\tilde{f}$   is non decreasing and $\tilde{f}(l(p,a)) \geq  \sum_{s\in S} \lambda_s\tilde{f}( l(p_s,a_s)).$ We obtain:
\begin{eqnarray*}
  \Phi(\alpha,f)(p)  & \geq & \alpha \sum_{s\in S} \lambda_s   g(p_s,a_s) + (1-\alpha)\sum_{s\in S} \lambda_s \tilde{f}( l(p_s,a_s)),\\
  \; & = &\sum_{s\in S} \lambda_s  \left( \alpha  g(p_s,a_s) + (1-\alpha) \tilde{f}( l(p_s,a_s))\right). \end{eqnarray*}
 \noindent This holds for every $(a_s)_{s\in S}$, so  $\Phi(\alpha,f)(p)   \geq  \sum_{s\in S} \lambda_s  \Phi(\alpha,f)(p_s).$ \hfill $\Box$

\subsubsection{Value of finite games and the recursive formula.}\label{subsub1}

\begin{lem} \label{lem4} For every state $p$ in $X$, the game $\Gamma_1(\delta_p)$ has a value which is: $$\dis v_1(p)= \max_{a\in A} \min_{b \in B} g(p,a,b)=\min_{b\in B} \max_{a \in A} g(p,a,b)= \Phi(1,0)(p).$$ \indent $v_1$ is concave and belongs to ${\cal D}$. $\tilde{v}_1$ is non decreasing and non expansive. \end{lem}

\noindent{\bf Proof:} Fix $p$ in $X$, and consider the game with normal form $(A,B, g(p,.,.))$.  By H3 and H4, we can apply    Sion's theorem (see e.g. Sorin 2002 p.156, thm A.7)  and obtain that this game has a value and both players have optimal strategies. By lemma  \ref{lem3,5}, we get that $v_1$ is concave, and by  H5,   we have   $v_1\in {\cal D}$. By claim \ref{cla2,8}, for every distribution $u$ the game $\Gamma_1(u)$ has a value which is precisely $\tilde{v}_1$. By concavity of $v_1$, $\tilde{v}_1$ is non decreasing. Since $v_1\in E_1$ and we use the Wasserstein distance, $\tilde{v}_1$ is non expansive. \hfill $\Box$

We will need to consider not only the $n$-stage games $\Gamma_n(u)$, but a larger family of games with initial distribution $u$.  

\begin{defi} \label{def9} Let   $\theta=\sum_{t\geq 1} \theta_t \delta_t$ be in $\Delta_f(\N^*)$, i.e. $\theta$ is a probability with finite support over positive integers. For $u$ in $\Delta_f(X)$, the game  $\Gamma_{[\theta]} (u)$ is the game with normal form $(\Sigma, {\cal T}, \gamma_{[\theta]}^{u})$, where: $$\gamma_{[\theta]} ^u (\sigma, \tau)=\E_{\P_{u, \sigma, \tau}} \left(  \sum_{t=1}^{\infty} \theta_t\; g (p_t,a_t,b_t)
\right).$$ \end{defi}
\noindent If $\theta=1/n\; \sum_{t=1}^n \delta_t$, $\Gamma_{[\theta]} (u)$ is nothing but $\Gamma_n(u)$. 
 $\Gamma_{[\theta]} (u)$ can be seen as the game where after the play, a stage   $t^*$ is selected according to $\theta$ and then only the payoff of stage $t^*$ matters. If $\theta=\sum_{t \geq 1} \theta_t \delta_t$, define $\theta^+$ as the law of $t^*-1$ given that $t^*\geq 2$. Define arbitrarily  $\theta^+=\theta$ if $\theta_1=1$, and otherwise we  have $\theta^+=\sum_{t \geq 1} \frac{\theta_{t+1}}{1-\theta_1}\delta_t$. We now write a recursive formula for the value of the games $\Gamma_{[\theta]} (u)$.  
 
\begin{pro} \label{pro2} For  $\theta=\sum_{t\geq 1} \theta_t \delta_t$ in $\Delta_f(\N^*)$ and $u$ in $\Delta_f(X)$, the game $\Gamma_{[\theta]} (u)$ has a value $\tilde{v}_{[\theta]}(u)$ such that:
\begin{eqnarray*} 
\forall p \in X, \;\;    v_{[\theta]}(p) & = &\Phi(\theta_1, {v}_{[\theta^+]})(p),\\
\; & = &\max_{a\in A}\; \theta_1 g(p,a) +(1-\theta_1)\tilde{v}_{[\theta^+]}(l(p,a)),\\
\; & =&\min_{b \in B}  \max_{a \in A}   \;  \theta_1 g(p,a,b) +(1-\theta_1) \tilde{v}_{[\theta^+]}(l(p,a)).
\end{eqnarray*}
 In $\Gamma_{[\theta]} (u)$, both players have optimal Markov strategies. $v_{[\theta]}$ is  concave and belongs to  ${\cal D}$. $\tilde{v}_{[\theta]}$ is non decreasing and non expansive. \end{pro}
  
  \noindent{\bf Proof:} by induction.  If $\theta=\delta_1$, lemma \ref{lem4} gives the result. 
  
  Fix now $n \geq2$, and assume that the proposition is true for every $\theta$ with support included in $\{1,...,n-1\}$. Fix a probability $\theta=\sum_{t=1}^n \theta_t \delta_t$, and notice that $\theta^+$ has a support included  in $\{1,...,n-1\}$. Fix also $p$ in $X$. 
  
  Consider the auxiliary zero-sum game $\Gamma'_{[\theta]} (p)$ with normal form $(A,B, f_{[\theta]}^p)$, where $f_{[\theta]}^p(a,b)= \theta_1 g(p,a,b) +(1-\theta_1)\tilde{v}_{[\theta^+]}(l(p,a)).$ We will apply Sion's theorem to this game. By H3, $A$ and $B$ are compact convex subsets of topological vector spaces. For every $a$, $(b \mapsto  f_{[\theta]}^p(a,b))$ is convex l.s.c. by H4. Consider now a convex combination $\lambda a + (1-\lambda) a'$ in $A$. By H6, we have $l(p, \lambda a +(1-\lambda)a')\succeq \lambda l(p,a) +(1-\lambda) l(p,a')$. By the induction hypothesis, $\tilde{v}_{[\theta^+]}$ is non decreasing, so:
  \begin{eqnarray*}
    \tilde{v}_{[\theta^+]} \left(l(p, \lambda a +(1-\lambda)a')\right) &\geq& \tilde{v}_{[\theta^+]} \left( \lambda l(p,a) +(1-\lambda) l(p,a')\right),\\
    \; & =&   \lambda \tilde{v}_{[\theta^+]}(l(p,a)) +(1-\lambda) \tilde{v}_{[\theta^+]}(l(p,a')).
    \end{eqnarray*}
  \noindent Since $g$ is concave in $a$, we obtain that $f_{[\theta]}^p$ is concave in $a$. Regarding continuity, by H6 and the induction hypothesis, $(a\longrightarrow \tilde{v}_{[\theta^+]}(l(p,a))$ is continuous. By H3, $(a\mapsto g(p,a,b))$ is u.s.c., so $(a\mapsto  f_{[\theta]}^p(a,b))$ is u.s.c.. By Sion's theorem $\Gamma'_{[\theta]} (p)$ has a value which is:  
  \begin{eqnarray*} 
  v'_{[\theta]}(p) &= & \max_{a \in A} \min_{b \in B}   \theta_1 g(p,a,b) +(1-\theta_1) \tilde{v}_{[\theta^+]}(l(p,a))\\
  \; & =&\min_{b \in B}  \max_{a \in A}    \theta_1 g(p,a,b) +(1-\theta_1) \tilde{v}_{[\theta^+]}(l(p,a)).
  \end{eqnarray*}
  
  Consider now the original game $\Gamma_{[\theta]} (p)$, and a strategy pair $(\sigma, \tau)$ in $\Sigma \times {\cal {T}}$. Write $a=\sigma_1(p)$, resp. $b=\tau_1(p)$, for the first action played by player 1, resp. player 2,  in $\Gamma_{[\theta]} (p)$. Denote by  $\sigma^+_{p,a,b}$ the continuation strategy issued from $\sigma$ after $(p,a,b)$ has occurred at stage 1. $\sigma^+_{p,a,b}$ belongs to $\Sigma$, and plays at stage $n$ after ($p_1$,$a_1$,$b_1$,...,$p_n$) what $\sigma$ plays at stage $n+1$ after ($p$,$a$,$b$,$p_1$,$a_1$,$b_1$,...,$p_n$). Similarly denote by  $\tau^+_{p,a,b}$ the continuation strategy issued from $\tau$ after $(p,a,b)$ has occurred at stage 1. It is easy to check that:
  $$\gamma_{[\theta]}^p(\sigma, \tau)= \theta_1 g(p,a,b) + (1-\theta_1) \gamma_{[\theta^+]}^{l(p,a)} (\sigma^+_{p,a,b}, \tau^+_{p,a,b}).$$
   \noindent  Consequently, in the game $\Gamma_{[\theta]} (p)$ player 1 can guarantee $ \max_{a \in A} \min_{b \in B}   \theta_1 g(p,a,b) +(1-\theta_1) \tilde{v}_{[\theta^+]}(l(p,a))$ by playing a Markov strategy. Similarly player 2 has a Markov strategy which guarantees $\min_{b \in B}  \max_{a \in A}   \theta_1 g(p,a,b) +(1-\theta_1) \tilde{v}_{[\theta^+]}(l(p,a)).$ Since the two quantities coincide, $\Gamma_{[\theta]} (p)$ has a value $v_{[\theta]}(p)=v'_{[\theta]}(p)$, and both players have Markov optimal strategies. 
  
  This implies that for every $u$ in $\Delta_f(X)$, the game  $\Gamma_{[\theta]} (u)$ has a value which is the affine extension $\tilde{v}_{[\theta]}(u)$, and both players have Markov optimal strategies in $\Gamma_{[\theta]} (u)$. ${v}_{[\theta]}=\Phi(\theta_1, {v}_{[\theta^+]})$ and ${v}_{[\theta^+]}$ is concave, so by lemma \ref{lem3,5} ${v}_{[\theta]}$ is concave, and $\tilde{v}_{[\theta]}$ is non decreasing. By H5 ${v}_{[\theta]}$ is in ${\cal D}$, so ${v}_{[\theta]}$ is 1-Lipschitz, and finally  $\tilde{v}_{[\theta]}$   is non expansive. \hfill $\Box$
  
  Among the games $\Gamma_{[\theta]} (u)$, the following family will play an important role and deserves a specific notation.
  
  \begin{defi}\label{def11}
  For $m\geq 0$ and $n\geq 1$,  $\Gamma_{m,n}(u)$  is the game with normal form $(\Sigma, {\cal T}, \gamma_{m,n}^{u})$.\end{defi}
  \noindent Recall that in definition \ref{def10,5}, we put  $\gamma_{m,n} ^{u} (\sigma, \tau)$ $=\E_{\P_{u, \sigma, \tau}} \left(  \frac{1}{n}\sum_{t=m+1}^{m+n}  \; g(p_t,a_t,b_t) \right)$ for each $(\sigma, \tau)$.  So $\Gamma_{m,n}(u)$ is nothing but  $\Gamma_{[\theta]} (u)$ with $\theta=1/n\; \sum_{t=m+1}^{m+n} \delta_t$. We can apply the previous proposition and denote the value of $\Gamma_{m,n}(u)$ by $v_{m,n}(u)$. $v_{0,n}$ is just the value of the $n$-stage game $v_n$, and for convenience we put $v_0=0$. We have for all $p$ in $X$, and positive $m$ and $n$:
  \begin{eqnarray*}
  v_n(p) & = &\Phi(1/n, v_{n-1}) =  \frac{1}{n}\max_{a \in A} \min_{b \in B} \left( g(p,a,b) +(n-1) \tilde{v}_{n-1}(l(p,a))\; \right),\\
\; &=&   \frac{1}{n}  \min_{b \in B} \max_{a \in A}\left( g(p,a,b) +(n-1) \tilde{v}_{n-1}(l(p,a))\; \right),\\
   v_{m,n}(p)& = &\Phi(0, v_{m-1,n}) =  \max_{a \in A}   \; \tilde{v}_{m-1, n}(l(p,a)).\end{eqnarray*}
   
  In $\Gamma_{m,n}(u)$, the players first play   $m$ stages to control the state, then they play $n$ stages   for payoffs.   Moreover, player 2 does not control the transitions, so he can play arbitrarily in the first $m$ stages.
  
  \begin{lem} \label{lem5} Fix  $n\geq 1$. There exists a Markov strategy $\tau=(\tau_t)_{t \geq 1}$ for player 2 such that  $\forall m \geq 0$, $\forall \tau'=(\tau'_t)_{t \geq 1}$ in ${\cal T}$:
  
  the  condition $(\forall l=1,...,n,\; \forall p \in X, \tau'_{m+l}(.,...,p)=\tau_l(p))$ implies that   for every $u$ in $\Delta_f(X)$, $\tau'$ is an optimal strategy  for player 2 in $\Gamma_{m,n}(u)$.\end{lem}
  
  \noindent{\bf Proof:} For each $t$ in $\{1,...,n\}$, define $\tau_t$ as the mapping which plays, if the current state is $p\in X$, an element $b$ in $B$ achieving the minimum in :
  $$\min_{b \in B}  \;\; \max_{a \in A}  \frac{1}{n-t+1} \left( g(p,a,b)+(n-t) \tilde{v}_{n-t}(l(p,a)) \right)= v_{n+1-t}(p).$$
 Using the previous recursive formula, one can show by induction that this construction of $\tau$ is appropriate. \hfill $\Box$
 
 \subsubsection{Player 2 can guarantee $v^*(u)$ in  $\Gamma(u)$.} \label{subsub2}
 
 We now consider the game with infinitely many stages $\Gamma(u)$. The following results are similar to propositions 7.7 and 7.8 in Renault, 2006.
 
 \begin{pro} \label{pro3} In $\Gamma(u)$, player 2 can guarantee with Markov strategies the quantity: $$\inf_{n \geq 1} \limsup_{T} \left( \frac{1}{T} \sum_{t =0}^{T-1} \tilde{v}_{nt,n}(u) \right).$$\end{pro}
 
 \noindent{\bf Proof:} Fix $n\geq 1$, and consider $\tau_1$,..., $\tau_n$ given by lemma \ref{lem5}.  Divide the set of stages $\N^*$ into consecutive blocks $B^1$, $B^2$,..., $B^m$,... of equal length $n$. By lemma \ref{lem5}, the cyclic strategy $\tau'=(\tau_1,..., \tau_n, \tau_1,...,\tau_n,\tau_1,...,\tau_n,....)$ is optimal for player 2 in the game $\Gamma_{nm,n}(u)$, for each $m\geq 0$. $\tau'$    is a Markov strategy for player 2, and for every strategy $\sigma$ of player 1 in $\Sigma$ we have:
 $$\forall m\geq 0,\; \E_{\P_{u,\sigma, \tau'}} \left( \frac{1}{n} \sum_{t \in B^{m+1}} g(p_t,a_t,b_t)\right) \leq \tilde{v}_{nm,n}(u),$$
$${\rm so }\;\; \forall M\geq 1,\; \E_{\P_{u,\sigma, \tau'}} \left( \frac{1}{nM} \sum_{t =1}^{nM} g(p_t,a_t,b_t)\right) \leq \frac{1}{M} \sum_{m=0}^{M-1} \tilde{v}_{nm,n}(u).$$
\noindent And since $n$ is fixed and payoffs are bounded, we obtain that player 2 can guarantee with   Markov strategies: $ \limsup_{M} \left( \frac{1}{M} \sum_{m=0}^{M-1} \tilde{v}_{nm,n}(u) \right)$ in $\Gamma(u)$. \hfill $\Box$

This proof also shows the following inequality.

\begin{lem}\label{lem6} $\forall u \in \Delta_f(X), \forall n \geq 1, \forall T\geq 1,\;\; \tilde{v}_{nT}(u)\leq \frac{1}{T} \sum_{t=0}^{T-1} \tilde{v}_{nt,n}(u).$\end{lem}

The following quantity will turn out to be the value of $\Gamma(u)$.

\begin{defi} \label{def12} For every $u$ in $\Delta_f(X)$, we define:
$$v^*(u)=\inf_{n \geq 1} \sup_{m \geq 0} \tilde{v}_{m,n}(u).$$\end{defi}

\noindent Since $v^*(u) \geq$ $\inf_{n \geq 1} \limsup_{T} \left( \frac{1}{T} \sum_{t =0}^{T-1} \tilde{v}_{nt,n}(u) \right)$, by proposition \ref{pro3} player 2 can also guarantee $v^*(u)$ with Markov strategies  in $\Gamma(u)$. By claim \ref{cla3}, $\overline{v}(u)=\min\{v \in \R, {\rm player \;2 \; guarantees \;}$ $ v\; {\rm in}\; \Gamma(u)\;\}$, so we now have the following inequality chain:
 $$\boxed{\underline{v}(u) \leq \liminf_N\tilde{v}_N(u)\leq  \limsup_N\tilde{v}_N(u)\leq \overline{v}(u)\leq v^*(u).}$$
 
 \subsubsection{Markov strategies for player 1.} \label{subsub3}

By H1 player 2 does not control the transition, so a Markov strategy $\sigma$ induces, together with the initial distribution $u$, a probability distribution $\P_{u,\sigma}$ over $(X \times A)^{\infty}$, i.e. over sequences of states and actions for player 1. For $u$ in $\Delta_f(X)$ and $\sigma_1:X \longrightarrow A$, we denote by $H(u, \sigma_1)$ the law of the state of stage 2 if the initial distribution is $u$ and player 1 plays at stage 1 according to $\sigma_1$. We denote by $G(u, \sigma_1)$ the payoff guaranteed by $\sigma_1$ at stage 1.   And we also define the continuation strategy $\sigma^+$.

\begin{notas} \label{not4} $\;$

\noindent $G(u, \sigma_1)= \sum_{p \in X} u(p) g(p, \sigma_1(p)) $, and $H(u, \sigma_1)= \sum_{p \in X} u(p) l(p, \sigma_1(p)) \in \Delta_f(X).$ 

\noindent If $\sigma=(\sigma_t)_{t \geq1}$  is in $\Sigma^{M}$, we write $\sigma^+$ for the Markov strategy $(\sigma_t)_{t \geq2}$. \end{notas}

We now concentrate on what player 1 can achieve in $\Gamma(u)$ and  completely forget player 2. We use similar notations as in definition \ref{def9}.

\begin{defi}\label{def12,5}  For   $\theta=\sum_{t\geq 1}\theta_t \delta_t$   in $\Delta_f(\N^*)$, $u$  in $\Delta_f(X)$ and $\sigma$ in  $\Sigma^{M}$, we put:
$$\gamma_{[\theta]}^{u}(\sigma)= \E_{\P_{u, \sigma}} \left(\sum_{t\geq 1} \theta_t \, g({p}_t, {a}_t)\right)= \sum_{t\geq 1} \theta_t \, \gamma_{[\delta_t]}^{u}(\sigma).$$\end{defi}

\noindent For simplicity, we write $\gamma_{[t]}$ instead of  $\gamma_{[\delta_t]}$ for the payoff induced at stage $t$. Clearly, we have for every $t\geq 2$:
$$\gamma_{[t]}^{u}(\sigma)= \gamma_{[{t-1}]}^{H(u,\sigma_1)}(\sigma^+).$$

\begin{lem} \label{lem7} $\gamma_{[\theta]}^{u}(\sigma)= \min_{\tau \in {\cal T}} \gamma_{[\theta]}^{u}(\sigma, \tau).$\end{lem}

The proof is easy, the minimum on the RHS being  achieved by a Markov strategy $\tau$  such that for every $t$ and $p_t$, $\tau_t(p_t)$ achieves the minimum in $b$ of the quantity $g(p_t, \sigma_t(p_t),b)$. As a corollary of lemma \ref{lem7}, we obtain that $\sup_{\sigma \in \Sigma^{M} } \gamma_{[\theta]}^{u}(\sigma)=$ $\sup_{\sigma  \in \Sigma^{M} }\inf_{\tau \in {\cal T}} \gamma_{[\theta]}^{u}(\sigma, \tau).$ By proposition \ref{pro2}, $\Gamma_{[\theta]}(u)$ has a value, so we get:

\begin{cor} \label{cor1} $\;$

For  every  $\theta$   in $\Delta_f(\N^*)$ and $u$  in $\Delta_f(X)$,   $\;\;\; \sup_{\sigma \in \Sigma^{M} } \gamma_{[\theta]}^{u}(\sigma)=\tilde{v}_{[\theta]}(u).$ \end{cor}

As in definition \ref{def11}, we now specify notations for a particular class of probabilities.

\begin{defi} \label{def13} For    $n\geq 1$, $m\geq 0$, $u$  in $\Delta_f(X)$ and $\sigma$ in  $\Sigma^{M}$, we put: 
$$\gamma_{m,n}^{u}(\sigma)= \E_{\P_{u, \sigma}} \frac{1}{n} \left(\sum_{t=1}^n  \, g({p}_{m+t}, {a}_{m+t})\right)= \frac{1}{n}\sum_{t=1}^n  \gamma_{[{m+t]}}^{u}(\sigma).$$\end{defi}

Finally, we consider  a situation where player 1 does not precisely know the length of the game.

\begin{defi} \label{def14} Fix $m\geq 0$, $n\geq 1$, and $u$ in $\Delta_f(X)$. We define:  $$w_{m,n}(u)= \sup_{\sigma \in \Sigma^M} \min \{\gamma_{m,t}^{u} (\sigma), t\in \{1,..,n\}\}.$$\end{defi}

The mappings $w_{m,n}$ will play an important role in the sequel, while applying corollary 3.8 of Renault, 2007. We will show   in corollary \ref{cor2} that they  are non expansive. To prove  corollary \ref{cor2} we will use the following  lemma \ref{lem8} and propositions  \ref{pro4}, \ref{pro5} and \ref{pro6}. We start  by defining  an auxiliary game\footnote{We proceed similarly as in section 6.2. of Renault, 2007. However the situation is more technical here, and it will not be possible  to apply a standard minmax theorem to the game ${\cal A}(m,n, u)$.}.

\begin{defi} \label{def15} For $m\geq 0$, $n\geq 1$, and $u$ in $\Delta_f(X)$, we define ${\cal A}(m,n, u)$ as the zero-sum   game with normal form  $(\Sigma^M$, $\Delta(\{1,...,n\}),f)$, where:
$$\forall \sigma \in \Sigma^M, \forall \theta \in \Delta(\{1,...,n\}),\; f(\sigma, \theta)= \sum_{t=1}^n \theta_t \gamma_{m,t}^
u( \sigma).$$ \end{defi}
We will prove later  that ${\cal A}(m,n, u)$ has a value which is $w_{m,n}(u)$ (see proposition \ref{pro6} below).  Notice that in general $f(., \theta)$ is not concave  in $\sigma$. However, we will show that it is {\it concave-like} in $\sigma$, i.e. that: $\forall \sigma', \sigma''$, $\forall \lambda \in [0,1]$, there exists $\sigma$ such that $\forall \theta$,  $f(\sigma, \theta)\geq \lambda f(\sigma', \theta)+(1-\lambda) f(\sigma'', \theta)$.  We start with a characterization for the partial order $\succeq$.

\begin{pro}\label{pro4} Let $u$ and $v$ be in $\Delta_f(X)$. Write $u=\sum_{p \in X} u(p)\delta_p$. The following conditions are equivalent:

$(i)$ $u \succeq v$, and

$(ii)$ For every $p$ such that $u(p)>0$, there exist $  S(p)\geq1$, $\lambda_1^p,...,\lambda_{S(p)}^p\geq 0$ and $q_1^p,...,q_{S(p)}^p$ in $X$ such that: $\sum_{s=1}^{S(p)} \lambda_s^p=1$, $\sum_{s=1}^{S(p)} \lambda_s^p q_s^p=p$, and $v=\sum_{p\in X} u(p) \left( \sum_{s=1}^{S(p)} \lambda_s^p \delta_{q_s^p} \right)$. \end{pro}

The proof can be easily deduced from a theorem of Loomis (see Meyer, 1966, T26 p.283), which deals with positive measures  on $X$. Notice that condition $(ii)$ can be seen as follows: $u$ is the law of some random variable $X_1$ with values in $X$, $v$ is the law of some random variable $X_2$ with values in $X$, and we have the martingale condition: $\E(X_2|X_1)=X_1$.\\

In general,  $\gamma_{[t]}^{u}(\sigma)$ is not a non decreasing function of $u$. However, we have the following property.

\begin{lem} \label{lem8} Let  $n\geq 1$, $u$ and $v$ be  in $\Delta_f(X)$ such that $u \succeq v$. For every $\sigma \in \Sigma^M$, there exists $\sigma' \in \Sigma^M$ such that:
$$\forall t \in \{1,...,n\},\; \gamma_{[t]}^{u}(\sigma') \geq  \gamma_{[t]}^{v}(\sigma).$$\end{lem}

\noindent {\bf Proof:} by induction on $n$.

If $n=1$, there exists $\sigma'$ such that $\gamma_{[1]}^u(\sigma')=\tilde{v}_1(u)$. Since $\tilde{v}_1$ is non decreasing, $\tilde{v}_1(u)\geq \tilde{v}_1(v)\geq \gamma_{[1]}^v(\sigma).$

Fix now $n\geq 1$, and assume that the lemma is proved for $n$. Fix $u$ and $v$ in $\Delta_f(X)$ with $u\succeq v$, and $\sigma$ in $\Sigma^M$. We have $u=\sum_{p \in X} u(p) \delta_p$, and by proposition \ref{pro4} it is possible to write $v=\sum_{p\in X} u(p) \left( \sum_{s=1}^{S(p)} \lambda_s^p \delta_{q_s^p} \right)$, with $\lambda_s^p\geq 0$, $\sum_{s=1}^{S(p)} \lambda_s^p=1$, and  $\sum_{s=1}^{S(p)} \lambda_s^p q_s^p=p$ for each $p$ such that $u(p)>0$. Define for every such $p$ and $s$, $a^p_s= \sigma_1(q_s^p)$. By the splitting hypothesis H7, for every $p$ one can find $a^p\in A$ such that:     $$\dis  l(p,a^p) \succeq \sum_{s=1} ^{S(p)} \lambda_s^p\, l(q_s^p ,a_s^p) \; \;\; \rm { and }\it  \;\; \;  g(p,a^p) \geq \sum_{s=1} ^{S(p)}\lambda_s^p \,  g(q_s^p,a_s^p).$$
\noindent We define what $\sigma'$ plays at stage 1 if the state is $p$ as:  $\sigma'_1(p)=a^p.$

We have: $\gamma_{[1]}^u(\sigma')=\sum_p u(p) g(p,a^p) \geq  \sum_p u(p) \sum_{s=1}^{S(p)} \lambda_s^p g(q_s^p,a_s^p)=\gamma_{[1]}^v(\sigma)$, and $H(u, \sigma'_1)= \sum_p u(p)
 l(p,a^p)\succeq \sum_p u(p) \sum_{s=1}^{S(p)} \lambda_s^p \,  l(q_s^p,a_s^p)=H(v, \sigma_1).$
 
 Since $H(u, \sigma'_1)\succeq  H(v, \sigma_1)$, we apply the induction hypothesis to the continuation strategy $\sigma+$. We obtain the existence of some Markov strategy $\tau=(\tau_t)_{t \geq 1}$ such that: $\forall t \in \{1,...,n\}$, $\gamma_{[t]}^{H(u, \sigma'_1)} (\tau) \geq \gamma_{[t]}^{H(v, \sigma_1)} (\sigma^+)$. Define $\sigma'_t=\tau_{t-1}$ for each $t\geq 2$. $\sigma'=( \sigma'_t)_{t \geq 1}$ is a Markov strategy for player 1, and satisfies: $\gamma_{[1]}^u(\sigma')\geq \gamma_{[1]}^v(\sigma)$, and for $t\in \{2,...,n+1\}$:
 
   $\gamma_{[t]}^u(\sigma')=\gamma_{[t-1]}^{H(u, \sigma'_1)}(\sigma'^+)$ $\geq \gamma_{[t-1]}^{H(v, \sigma_1)} (\sigma^+)$ $= \gamma_{[t]}^v(\sigma)$. \hfill $\Box$
   
 \vspace{0,5cm}
 
  Lemma \ref{lem8} is now improved.
    
  \begin{pro} \label{pro5}   Let  $n\geq 1$, $\lambda \in [0,1]$, $u$, $u'$  and $u''$    in $\Delta_f(X)$ be such that $u \succeq \lambda u'+ (1-\lambda) u''$. For every $\sigma'$ and $\sigma''$ in $\Sigma^M$, there exists $\sigma\in \Sigma^M$ such that:
$$\forall t \in \{1,...,n\},\; \gamma_{[t]}^{u}(\sigma) \geq  \lambda \gamma_{[t]}^{u'}(\sigma')+ (1-\lambda) \gamma_{[t]}^{u''}(\sigma'').$$\end{pro}

\noindent {\bf Proof:} by induction on $n$.

If $n=1$, there exists $\sigma$ such that  $\gamma_{[1]}^{u}(\sigma)$ $ = \tilde{v}_1(u) \geq  \tilde{v}_1(\lambda u'+(1-\lambda) u'') $ $=\lambda \tilde{v}_1(u') +(1-\lambda) \tilde{v}_1(u'') \geq \lambda \gamma_{[1]}^{u'} (\sigma') +(1-\lambda) \gamma_{[1]}^{u''}(\sigma'').$

Assume the proposition is proved for some $n\geq1$, and fix $u$, $u'$, $u''$, $\lambda$, $\sigma'$ and $\sigma''$  with $u \succeq \lambda u'+ (1-\lambda) u''$. Write $v= \lambda u'+ (1-\lambda) u''$. By lemma \ref{lem8} it is enough to find $\sigma$ in $\Sigma^M$ such that: $\forall t\in \{1,...,n+1\}$, $\gamma_{[t]}^{v}(\sigma) \geq  \lambda \gamma_{[t]}^{u'}(\sigma')+ (1-\lambda) \gamma_{[t]}^{u''}(\sigma'').$

We have $v=\sum_p \left(\lambda u'(p) +(1-\lambda) u''(p)\right) \delta_p$, and $v(p)=\lambda u'(p) +(1-\lambda) u''(p)$ for each $p$. For every $p$ such that $v(p)>0$, we define:
$$\sigma_1(p)= \frac{\lambda u'(p)}{v(p)} \sigma'_1(p) +\frac{(1-\lambda)u''(p)}{v(p)} \sigma''_1(p).$$
\noindent  $\sigma_1(p)$ belongs to $A$ by convexity. Now,
\begin{eqnarray*}
\lambda \gamma_{[1]}^{u'} (\sigma') +(1-\lambda) \gamma_{[1]}^{u''}(\sigma'') & =& \lambda \sum_{p} u'(p) g(p,\sigma_1'(p)) +(1-\lambda) \sum_{p} u''(p) g(p, \sigma_1''(p)) \\
\; & = & \sum_{p} v(p) \left( \frac{\lambda u'(p)}{v(p)} g(p,  \sigma'_1(p)) +\frac{(1-\lambda)u''(p)}{v(p)}g(p,  \sigma''_1(p)) \right) \\
\; & \leq & \sum_{p} v(p)  \; g(p,    \sigma_1(p) )=\gamma_{[1]}^v(\sigma),
 \end{eqnarray*}
\noindent where the inequality comes from the concavity of $g$ in the variable $a$ (see H4). Proceeding   in the same way with distributions on the second state, we obtain via the concavity of $l$ in $a$ (see H6):

\begin{eqnarray*}
\lambda H({u'}, \sigma'_1) +(1-\lambda)H({u''},\sigma''_1) & =& \lambda \sum_{p} u'(p) l(p,\sigma_1'(p)) +(1-\lambda) \sum_{p} u''(p) l(p, \sigma_1''(p)) \\
\; & = & \sum_{p} v(p) \left( \frac{\lambda u'(p)}{v(p)} l(p,  \sigma'_1(p)) +\frac{(1-\lambda)u''(p)}{v(p)}l(p,  \sigma''_1(p)) \right) \\
\; & \preceq & \sum_{p} v(p)  \; l(p,     \sigma_1(p) )=H(v,\sigma)\\
\end{eqnarray*}

\noindent Consequently, we have   $H(v,\sigma)\succeq\lambda H({u'}, \sigma'_1) +(1-\lambda)H({u''},\sigma''_1)$, and by the induction hypothesis  there exists $\sigma^+$ in $\Sigma^M$ such that $\forall t \in \{1,...,n\},\; \gamma_{[t]}^{H(v, \sigma_1)}(\sigma^+) \geq  \lambda \gamma_{[t]}^{H(u', \sigma'_1)}(\sigma'^+)+ (1-\lambda) \gamma_{[t]}^{H(u'', \sigma_1'')}(\sigma''^+).$ We  naturally define $\sigma=(\sigma_1, \sigma^+)$, and we have, for $t\in \{2,...,n+1\}$: $\gamma_{[t]}^v(\sigma)=\gamma_{[t-1]}^{H(v, \sigma_1)}(\sigma^+)$ $\geq \lambda \gamma_{[t-1]}^{H(u', \sigma'_1)} (\sigma'^+) + (1-\lambda) \gamma_{[t-1]}^{H(u'', \sigma''_1)} (\sigma''^+)$ $= \lambda \gamma_{[t]}^{u'}(\sigma')+ (1-\lambda) \gamma_{[t]}^{u''}(\sigma'').$
 \hfill $\Box$
 
   \begin{pro} \label{pro6}  For every $m\geq 0$, $n\geq 1$ and $u$ in $\Delta_f(X)$, the game ${\cal A}(m,n, u)$ has a value which is $w_{m,n}(u)$.\end{pro}
   
   \noindent{\bf Proof:} Recall that the payoff function in  ${\cal A}(m,n, u)$ is:   $f(\sigma, \theta)= \sum_{t=1}^n \theta_t \gamma_{m,t}^
u( \sigma)$ for every $\sigma$ in $\Sigma^M$ and $\theta$ in $\Delta(\{1,...,n\})$. 

$\Delta(\{1,...,n\})$ is convex and compact, and $f$ is affine continuous in $\theta$. We now show that  $f$ is  concave-like in $\sigma$. Let $\sigma', \sigma''$ be in $\Sigma^M$, and $\lambda \in [0,1]$. By the previous proposition, there exists $\sigma$ in $\Sigma^M$ such that: $\forall t \in \{1,...,m+n\}$, $\gamma_{[t]}^u(\sigma) \geq \lambda \gamma_{[t]}^u(\sigma') + (1-\lambda) \gamma_{[t]}^u(\sigma'').$ So $f(\sigma, \theta)=\sum_{t=1}^n \frac{\theta_t}{t} \sum_{t'=1}^t \gamma_{[m+t']}^u(\sigma)$ $\geq$  $\sum_{t=1}^n \frac{\theta_t}{t} \sum_{t'=1}^t \left( \lambda \gamma_{[m+t']}^u(\sigma') + (1-\lambda) \gamma_{[m+t']}^u(\sigma'')\right)$ $= \lambda f(\sigma', \theta) +(1-\lambda) f(\sigma'', \theta).$ By a theorem of Fan (1953, see proposition A.13 p.160 in Sorin 2002), ${\cal A}(m,n, u)$ has a value which is:
$\sup_{\sigma \in \Sigma^M} \inf_{\theta \in \Delta(\{1,...,n\})} f(\sigma, \theta)$ $= w_{m,n}(u).$ \hfill $\Box$

\begin{nota} \label{nota3} For $\theta=\sum_{t=1}^n \theta_t \delta_t$ $\in \Delta(\{1,...,n\})$, and $m\geq0$, we define $\theta^{m,n}$ in $\Delta(\{1,...,m+n\})$ by: 

\centerline{$\theta^{m,n}_s=0$ if $s\leq m$, and  $\theta^{m,n}_s=\sum_{t=s-m}^n \frac{\theta_t}{t}$ if  $m<s\leq n+m$.} \end{nota}

\begin{cor} \label{cor2}For every $m\geq0$, $n\geq 1$ and $u$ in $\Delta_f(X)$,
$$w_{m,n}(u) =\inf_{\theta\in \Delta(\{1,...,n\})} \tilde{v}_{[\theta^{m,n}]}(u).$$
The mapping $w_{m,n}$ is non decreasing and non expansive.\end{cor}

\noindent {\bf Proof:} Fix $m$, $n$ and $u$. $w_{m,n}(u)$ is the value of ${\cal A}(m,n, u)$, so we have:
$ w_{m,n}(u) =\inf_{\theta\in \Delta(\{1,...,n\})}\sup_{\sigma \in \Sigma^M}  f(\sigma, \theta)$.  But $f(\sigma, \theta)=\sum_{t=1}^n \theta_t \gamma_{m,t}^{u} (\sigma)= \gamma_{[\theta^{m,n}]}^{u}(\sigma).$ So we obtain: $ w_{m,n}(u) =\inf_{\theta\in \Delta(\{1,....,n\})}\sup_{\sigma \in \Sigma^M}  \gamma_{[\theta^{m,n}]}^{u}(\sigma).$ Corollary \ref{cor1} gives $w_{m,n}(u) =\inf_{\theta\in \Delta(\{1,....,n\})} \tilde{v}_{[\theta^{m,n}]}(u).$ For each $\theta$, $\tilde{v}_{[\theta^{m,n}]}$ is  non decreasing and non expansive by proposition \ref{pro2}, hence the result. \hfill $\Box$

\vspace{1cm}

\subsubsection{A dynamic programming problem.} \label{subsub4}

We will   conclude the proof of theorem \ref{thm3} and show that the uniform value exists.  Fix the initial distribution $u$, our stochastic game is defined by $\Gamma(u)=(X,A,B,g,l,u).$ We now define an auxiliary MDP    as follows.\\

\begin{defi} \label{def16} The MDP  $\Psi(z_0)$ is defined as $(Z,F,r,z_0)$, where: 
\begin{quote}
 $\bullet$ $Z=\Delta_f(X) \times [0,1]$ is the set of states of the MDP,
 
  $\bullet$  $z_0=(u,0)$ is the initial state in $Z$,
    
  $\bullet$     $r$ is the payoff function from $Z$ to $[0,1]$ defined by $r(u,y)=y$ for each $(u,y)$ in $Z$, 
       
  $\bullet$     and   the transition function $F$ is the correspondence from $Z$ to $Z$ such that: $\forall z=(u,y)\in Z$, 
\begin{eqnarray*}
  F(z) & =& \left\{\left( \sum_{p \in X} u(p) l(p,a(p)), \sum_{p \in X} u(p) g(p,a(p))\right), \; \forall p \; a(p) \in A \right\},\\
  \; & = &\{ (H(u,f), G(u,f)), f { \; is \;  a \; mapping \; from\;  } X { \; to \; } A\}.
  \end{eqnarray*}\end{quote}\end{defi}
 
 Notice that $F(u,y)$ does not depend on $y$, hence the  value functions will not depend on $y$. $F$ has non empty values. Even with strong assumptions on $l$ and $g$, $F$ may not have a compact graph, because in the definition of $F(z)$ we have a unique $a(p)$ for each $p$. So even if $q$ is close to $p$, the image by $F$ of $(1/2 \delta_p + 1/2 \delta_q,0)$ may be quite larger than $F(\delta_p,0)$.
 
 As in   Renault, 2007, we denote by $S(z_0)=\{s=(z_1,...,z_t,...)\in Z^{\infty}, \forall t\geq 1, z_t\in F(z_{t-1})\}$  the set of plays at $z_0$. The next proposition shows the strong links between the stochastic game $\Gamma(u)$ and the MDP $\Psi(z_0)$.
 
 \begin{pro} \label{pro7}
 a) For every Markov strategy $\sigma$ in $\Sigma^M$, there exists a play $s =(z_1,...,z_t,...)\in S(z_0)$ such that $\forall t\geq 1,\; \gamma_{[t]}^u(\sigma)=r(z_t).$
 
 b) Reciprocally, for every play $s =(z_1,...,z_t,...)\in S(z_0)$, there exists a Markov strategy $\sigma$ in $\Sigma^M$ such that: $\forall t\geq 1,\; \gamma_{[t]}^u(\sigma)=r(z_t).$
 \end{pro}
 
 \noindent{\bf Proof:} a) Take a Markov strategy $\sigma=(\sigma_t)_{t\geq 1}$ in $\Sigma^M$. Put   $u_1=u$ and $y_0=0$,  so that  $z_0=(u_1,y_0)$. Define by induction, for every $t\geq 1$, $u_{t+1}=H(u_t,\sigma_t)$, $y_t=G(u_t, \sigma_t)$, and $z_{t}=(u_{t+1},y_t)\in F(z_{t-1})$. $s=(z_t)_{t\geq 1}$ is a play at $z_0$. 
 
 $\gamma_{[1]}^u(\sigma)=G(u, \sigma_1)=y_1=r(z_1)$, and for $t\geq 2$,  $\gamma_{[t]}^u(\sigma)=\gamma_{[t-1]}^{H(u,\sigma_1)}(\sigma^+)$ $=$ $\gamma_{[t-1]}^{u_2}(\sigma^+)$ $=\gamma_{[t-2]}^{u_3}(\sigma^{++})$ $=...$ $=$ $\gamma_{[1]}^{u_t}((\sigma_{t'})_{t'\geq t})$ $=G(u_t,\sigma_t)=y_t=r(z_t)$. 
 
 b) Take a play $s=(z_t)_{t \geq 1}$ at $z_0$. Write for each $t\geq 0$, $z_t=(u_{t+1},y_t)\in \Delta_f(X) \times [0,1]$. For every $t\geq 1$, there exists a mapping $f_t$ from $X$ to $A$ which defines $z_t=(u_{t+1},y_t)$ in terms of $z_{t-1}$, i.e. that $u_{t+1}=H(u_t,f_t)$, and $y_t=G(u_t, f_t)$.  Simply define $\sigma=(f_t)_{t \geq 1}$.  As in point a), one can check that $\gamma_{[1]}^u(\sigma)=r(z_t)$ for each positive $t$. \hfill $\Box$

 \vspace{0,5cm}

 For any $m\geq 0$,  $n\geq 1$, and $s=(z_t)_{t\geq 1}$   in $Z^{\infty}$, we put as in definitions   3.1. 3.2.  of Renault, 2007: 

 $ \gamma_{m,n}(s) ={1 \over n} \sum_{t=1}^n r(z_{m+t}),$ 
 
  $\nu_{m,n}(s) = \min \{\gamma_{m,t}(s), t \in \{1,...,n\}\},$
  
   $v_{m,n}(z_0)= \sup_{s\in S(z_0)} \gamma_{m,n}(s),$ and 
   
   $ w_{m,n}(z_0)= \sup_{s\in S(z_0)} \nu_{m,n}(s).$

  By proposition \ref{pro7} and corollary \ref{cor1}, it is easy to obtain for every $m$ and $n$,  the equality of the values in the stochastic game $\Gamma(u)$ and in the MDP $\Psi(z_0)$: $\tilde{v}_{m,n}(u)=v_{m,n}(z_0)$. As  a  consequence,   we also have $v^*(u)=\inf_{n \geq 1}\sup_{m \geq 0} \tilde{v}_{m,n}(u)=v^*(z_0)$ (see definitions 3.6 of Renault, 2007 and \ref{def12} here). Similarly, we have $w_{m,n}(z_0)=$ $ \sup_{\sigma \in \Sigma^M} $ $\min \{\gamma_{m,t}^{u} (\sigma), t\in \{1,..,n\}\}=w_{m,n}(u)$ (see definition \ref{def14}).
 
 Define now $\underline{v}^M(u)$ as the maximal quantity that can be guaranteed by player 1 in $\Gamma(u)$ with Markov strategies: $$\underline{v}^M(u) = \sup_{\sigma\in \Sigma^M} \liminf_n \left( \inf_{\tau \in {\cal T} }\gamma_n^u( \sigma, \tau)\right)= \sup_{\sigma\in \Sigma^M} \liminf_n  \gamma_{0,n}^u( \sigma).$$
  Recall   that the  lower   value of the MDP is defined by: 
  
  $\underline{v}(z_0)=  \sup_{(z_t)_{t\geq 1} \in S(z_0)}  \left( \liminf_n  \frac{1}{n} \sum_{t=1}^n r(z_{t}) \right).$   
Again, proposition \ref{pro7} gives the equality $\underline{v}^M(u)=\underline{v}(z_0)$, so that we have the following relations:\\

   \begin{tabular}  {|c|}
   \cline{1-1}
$\underline{v}(z_0)=  \underline{v}^M(u)\leq \underline{v}(u)  \leq \liminf_N\tilde{v}_N(u)=\liminf_N {v}_N(z_0)$\\
 $ \leq  \limsup_N {v}_N(z_0)=\limsup_N\tilde{v}_N(u)  \leq \overline{v}(u)\leq v^*(u)=v^*(z_0). $\\
   \cline{1-1}
  \end{tabular}  
  
  \vspace{0,5cm}

We can now conclude.  We use the Wasserstein distance on $\Delta_f(X)$,  so $Z$ naturally is a precompact metric space. For every $m$ and $n$, by corollary \ref{cor2},  $(u \mapsto w_{m,n}(u))$ is a non expansive mapping from $\Delta_f(X)$ to $[0,1]$. This implies that $(z_0\mapsto w_{m,n}(z_0))$ is a non expansive mapping from $Z$ to $[0,1]$. As a consequence the family $(w_{m,n})_{m\geq 0, n \geq1}$ is uniformly continuous. By corollary 3.8  of Renault, 2007, we obtain that the MDP $\Psi(z_0)$ has a uniform value which is: $$v^*(z_0)=  \underline{v}(z_0) =  \lim_N {v}_N(z_0)= \sup_{m\geq 0} \inf_{n\geq 1} w_{m,n}(z_0)= \sup_{m\geq 0} \inf_{n\geq 1} v_{m,n}(z_0).$$ And the convergence from $(v_n)$ to $v^*$ is uniform. 
Back to our stochastic game $\Gamma(u)$, we obtain that $(v_n)_n$ uniformly converges to $v^*$, and 
 $$ v^*(u)=\underline{v}^M(u)=\underline{v}(u)= \lim_N\tilde{v}_N(u) = \overline{v}(u).$$ 
 So $\underline{v}(u)=  \overline{v}(u)$, which implies that $\Gamma(u)$ has a uniform value. Moreover,  $v^*(u)=$ $\sup_{m\geq 0} \inf_{n\geq 1} w_{m,n}(u)$ $ = \sup_{m\geq 0} \inf_{n\geq 1} \tilde{v}_{m,n}(u)$, and using definition 3.6 of Renault, 2007, 
 $v^*(u)=$ $ \inf_{n\geq 1} $ $ \sup_{m\geq 0} w_{m,n}(u)$ $=  \inf_{n\geq 1} \sup_{m\geq 0} \tilde{v}_{m,n}(u)$.
 This concludes the proof of theorem \ref{thm3}.
 
\centerline{$\Box$ \hspace{0,2cm} $\Box$ \hspace{0,2cm} $\Box$ \hspace{0,2cm} $\Box$}

\subsection{Comments.}\label{sub5,1}

\begin{rem} \label{rem5} Player 2 has 0-optimal strategies.

Under the same hypotheses H1,..., H6, it is possible to slightly modify the proof of Proposition \ref{pro3} and obtain that player 2 has a Markov strategy $\tau$ which is 0-optimal in $\Gamma(u)$, i.e. such that: $\forall \varepsilon>0,   \exists N_0, \forall N\geq N_0, \forall \sigma \in \Sigma, $   $  \gamma_N^u(\sigma, \tau) \leq v+ \varepsilon.$ 

Divide the set of stages into consecutive blocks $B^1$,..., $B^m$,..., such that $B^m$ has cardinal $m$  for each $m$. By lemma \ref{lem5}, there exists a Markov strategy $\tau=(\tau_t)_{t \geq 1}$ with the property that for every $m>0$, $\tau$ plays optimally within $B^m$, in the sense that  $\tau$ is an optimal strategy for player 2 in $\Gamma_{m(m-1)/2,m}(u)$.  For every strategy $\sigma$ of player 1, we have: 
$$\forall m\geq 1,\; \E_{\P_{u,\sigma, \tau}} \left( \frac{1}{m(m+1)/2} \sum_{t =1}^{\max(B^m)} g(p_t,a_t,b_t)\right) \leq \frac{1}{m} \sum_{i=1}^{m} \tilde{v}_{i(i-1)/2,i}(u).$$
\noindent We have seen in subsubsection \ref{subsub4}  that the values in the stochastic game $\Gamma(u)$ are the values of the  MDP $\Psi(z_0)$, so we can apply lemma 3.4 of Renault, 2007:  $\forall i\geq 1, \forall k \geq 1,$ $ \tilde{v}_{i(i-1)/2,i}(u)\leq \sup_{l \geq 0} w_{l,k}(u) +\frac{k-1}{i}.$ 

Fix now $\varepsilon>0$. One can find $k$ such that $ \sup_{l \geq 0} w_{l,k}(u)\leq v^*(u) + \varepsilon$. Since $\frac{k-1}{i} \longrightarrow_{i \to \infty}0$, one can find $m_0$ such that for every $m\geq m_0$, $$\E_{\P_{u,\sigma, \tau}} \left( \frac{1}{m(m+1)/2} \sum_{t =1}^{\max(B^m)} g(p_t,a_t,b_t)\right) \leq \frac{1}{m} \sum_{i=1}^{m} \left( v^*(u)+ \varepsilon + \frac{k-1}{i} \right) \leq v^*(u)+ 2 \varepsilon.$$ Looking at the size of the blocks,  one can show that $\tau$ is 0-optimal for player 2.  \hfill $\Box$
\end{rem}

\vspace{0,5cm}

\begin{rem} \label{rem5,3} A simple hypothesis implying H5.

Recall that hypothesis H5  requires the existence of some  subset  ${\cal D}$ of $E_1$  which contains  $\Phi(1,0)$, and is stable under any $ \Phi(\alpha,.)$. The following hypothesis is   stated in terms of the mappings $g$ and $l$. The distance $d$ on $X$ is extended to $\Delta(X)$ by the Wasserstein distance. \\

 \noindent {\underline{Hypothesis H5'}}: $\forall p \in X, \forall a \in A, \forall p' \in X, \exists a'\in A$ such that: 
 
 $d(l(p,a), l(p',a'))\leq d(p,p')$, and $\inf_{b \in B} g(p',a',b)\geq \inf_{b \in B} g(p,a,b) -d(p,p')$.\\
 
 It is easy to check that H5' implies that $E_1$ itself is stable under any  $ \Phi(\alpha,.)$ so H5' implies H5.  Consequently, the conclusions of theorem \ref{thm3} are true  if H1,H2,H3,H4,H5',H6, H7 hold.\end{rem}

 \vspace{2cm}

\section{Repeated games with an informed controller}\label{sec45}

\subsection{Model} \label{subsub37} We first consider a general model of zero-sum repeated game. We have: 
 \begin{quote}

 $\bullet$ five  non empty finite  sets:  a set of states or parameters  $K$,    a set $I$ of actions for player 1, a set $J$ of actions   for player 2,  a set $C$ of signals for player 1 and a set $D$ of signals for player 2.
 
 $\bullet$ an  initial  distribution $\pi \in \Delta(K \times C \times D)$,  
  
$\bullet$  a mapping $g$ from $K \times I \times J$ to $[0,1]$, called the payoff function of player 1, and
 
 $\bullet$ a mapping $q$ from $K \times I \times J$ to   $\Delta(K \times C \times D)$, called the transition function.  

 \end{quote}
 
 \noindent    Initially, $(k_1,c_1, d_1)$ is selected according to $\pi$,  player 1 learns $c_1$ and  player 2 learns $d_1$. Then simultaneously player 1 chooses $i_1$ in $I$ and player 2 chooses $j_1$ in $J$. The payoff for player 1 is $g(k_1,i_1,j_1)$,  then $(k_2,c_2,d_2)$ is selected according to $q(k_1,i_1,j_1)$, etc...  At any stage $t\geq2$, $(k_t, c_t, d_t)$ is selected    according to $q(k_{t-1},i_{t-1},j_{t-1})$,  player 1 learns $c_t$ and  player 2 learns $d_t$. Simultaneously, player 1 chooses $i_t$ in $I$ and player 2 chooses $j_t$ in $J$. The stage payoffs are $g(k_t,i_t,j_t)$ for player 1 and the opposite for player 2, and the play proceeds to stage $t+1$. \\
 
From now on we fix $\Gamma=(K,I,J,C,D,g,q)$, and for every $\pi$ in $\Delta(K \times C \times D)$ we denote by $\Gamma(\pi)=(K,I,J,C,D,\pi, g,q)$ the corresponding repeated game induced by $\pi$. For the moment we make no assumption on $\Gamma$, so we have a general model including stochastic games, repeated games with incomplete information and imperfect monitoring (signals). We   start with elementary definitions and notations. 

Players are allowed to select their actions randomly. A (behavior) strategy for player 1 is a sequence $\sigma=(\sigma_t)_{t \geq 1}$, where for each $t$, $\sigma_t$ is a mapping from $(C \times I)^{t-1}\times C$ to $\Delta(I)$, with the interpretation that $\sigma_t(c_1,i_1,...,c_{t-1},i_{t-1},c_t)$ is the lottery on  actions used by player 1 at stage $n$ after $(c_1,i_1,...,c_{t-1},i_{t-1},c_t)$. $\sigma_1$ is just a mapping from $C$ to $\Delta(I)$ giving the first action played by player 1 depending on his  initial signal. Similarly, a strategy for player 2 is a  sequence $\tau=(\tau_t)_{t \geq 1}$, where for each $t$, $\tau_t$ is a mapping from $(D \times J)^{t-1}\times D$ to $\Delta(J)$. We denote by $\Sigma$ and ${\cal T}$ the sets of strategies of player 1 and player 2, respectively.

It is standard that a pair of strategies $(\sigma, \tau)$ induces a   probability $\P_{\pi,\sigma, \tau} $ on the set of plays $\Omega=(K\times C \times D\times I \times J)^{\infty}$, endowed with the $\sigma$-algebra generated by the cylinders.  

\begin{defi}\label{def17}   The payoff for player 1
 induced by $(\sigma,\tau)$ at the first $N$ stages   is denoted by: 
$$\gamma_N^\pi (\sigma, \tau)=\E_{\P_{\pi, \sigma, \tau}} \left( {1 \over N}\sum_{n=1}^N g (k_n,i_n,j_n)
\right).$$   \end{defi}

  For $\pi$ in $\Delta(K \times C \times D)$ and $N\geq 1$,   the  $N$-stage  game game  $\Gamma_N(\pi)$ is the zero-sum game with normal form $(\Sigma, {\cal T}, \gamma_N^{\pi})$. By Kuhn 's theorem, $\Gamma_N(\pi)$ can be seen as the mixed extension  of a  finite game, so it has a value $v_n(\pi)$. 
  
  The following definitions are similar to those of subsection \ref{sub29} or section 2 of Renault, 2007.

 \begin{defi}\label{de18} Let $\pi$ be in $\Delta(K \times C \times D)$.
 
$$\mathnormal{ The \; liminf \; value\;  of } \; \Gamma (\pi) \; \mathnormal{ is:} \; v^-(\pi)= \liminf_n v_n(\pi).$$ 
$$\mathnormal{  The \; limsup\; value\;  of } \; \Gamma (\pi) \; \mathnormal{ is:} \; v^+(\pi)= \limsup_n v_n(\pi).$$

The  lower (or maxmin)  value of $\Gamma(\pi)$ is: $$\underline{v}(\pi) = \sup_{\sigma \in \Sigma} \liminf_n \left( \inf_{\tau \in {\cal T} }\gamma_n^{\pi}( \sigma, \tau)\right).$$

 The  upper (or minmax)  value of $\Gamma(u)$ is:
  $$\; \overline{v}(\pi) = \inf_{\tau \in {\cal T} }\limsup_n \left(\sup_{\sigma\in \Sigma}  \gamma_N^\pi( \sigma, \tau)\right).$$
 
$$\mathnormal{We\;  have}\;    \underline{v}(\pi)\leq v^-(\pi) \leq v^+(\pi) \leq \overline{v}(\pi).$$  $\Gamma(\pi)$ is said to have a uniform value if and only if $\underline{v}(\pi)=\overline{v}(\pi)$, and in this case the uniform value is $\underline{v}(\pi)=\overline{v}(\pi)$. \end{defi}

An equivalent definition of the uniform value is as follows. Given a real number $v$,  we say that player 1 can guarantee $v$ in $\Gamma(\pi)$  if: $\forall \varepsilon>0, \exists \sigma \in \Sigma, \exists N_0, \forall N\geq N_0, \forall \tau\in {\cal T}, \gamma_N^\pi(\sigma, \tau) \geq v- \varepsilon.$ Player 2 can guarantee $v$ in $\Gamma(\pi)$ if: $\forall \varepsilon>0, \exists \tau\in {\cal T}, \exists N_0, \forall N\geq N_0, \forall \sigma \in \Sigma, \gamma_N^\pi(\sigma, \tau) \leq v+ \varepsilon.$ If player 1 can guarantee $v$ and player 2 can  guarantee $w$ then clearly $w\geq v$. We also have, exactly as in claim \ref{cla3}: 
\begin{cla} \label{cla5} 

 $\underline{v}(\pi)=\max\{v \in \R, {\rm player \; 1\; can \; guarantee\; } v\; {\rm in} \;\Gamma(\pi)\;  \},$ 

\hspace{1,7cm}$\overline{v}(\pi)=\min\{v \in \R, {\rm player \;2 \; can \; guarantee \;} v\; {\rm in}\; \Gamma(\pi)\;\}.$

A real number  $v$ can be guaranteed by both players if and only if  $v$ is   the uniform value of $\Gamma(\pi)$.
\end{cla}
 
We now consider hypotheses on $q$ and $\pi$.\\

\noindent {\bf Hypothesis HA}: Player 1 is    informed of everything, i.e. at any stage $t\geq 1$, 
  he can deduce from his  signal $c_t$: the state $k_t$, player 2's signal $d_t$, and if $t\geq 2$, he can also deduce from $c_t$ the action $j_{t-1}$ previously played by player 2.\\
  
 \noindent {\bf Hypothesis HB}: Player 1 fully controls the transition, i.e. $q(k,i,j)$ does not depend on $j$ for each $(k,i)$ in $K \times I$.\\
 
HA and HB are very strong hypotheses, and they are incompatible as soon as $J$ has several elements. We will use  weaker hypotheses.\\

\noindent {\bf Hypothesis HA'}: Player 1 is    informed, in the sense that he can always deduce the state and player 2's signal from his own signal. Formally, there exists two mappings $\hat{k}:C \longrightarrow K$ and $\hat{d}:C \longrightarrow D$ such that, if  $E$ denotes $\{(k,c,d)\in K \times C \times D, \hat{k}(c)=k\; \mbox{ and } \hat{d}(c)=d\}$, then:$$\pi(E)=1, \; {\rm{and}} \; \; q(k,i,j)(E)=1, \; \forall (k,i,j)\in K \times I \times J.$$
 
 Notice that HA' does not mean that player 1 knows everything. Since we did not include the signals in the move, knowing the signal $d_t$ of player 2 at some stage $t$ does not imply knowing the action $j_t$ by player 2. However, not knowing this action will not be a problem for player 1 because we will also assume that player 2 does not really influence the transitions.\\
 
\noindent {\bf Hypothesis HB'}: Player 1 controls the transition, in the sense that the marginal of the transition $q$ on $K \times D$ does not depend on player 2's action. For $k$ in $K$, $i$ in $I$ and $j$ in $J$, we denote by $\bar{q}(k,i)$ the marginal of $q(k,i,j)$ on $K\times D$.\\

Assume that HA' and HB' hold. The couple (new state, signal of player 2) is selected according to a distribution depending on the current state and player 1's action, but not depending on player 2's action.  Player 2 may influence the distribution of player 1's signal, but still  player 1 will    be able to deduce the state and player 2's new information on the state. So essentially player 2 can influence player 1's knowledge about player 2's action. But this information is not relevant because it does not affect player 2's belief on the future  states.

 \begin{thm} \label{thm5} Under the hypotheses HA' and HB', the repeated game $\Gamma(\pi)$ has a uniform value. \end{thm}
  
  The next subsection is devoted to the proof of theorem \ref{thm5}. See subsection \ref{subcomcon} for other comments on hypotheses, applications and open questions. 
 
 \subsection{Proof of theorem \ref{thm5}}
 
 We   assume in this subsection  that HA' and HB' are satisfied.  Keeping fixed all other quantities, increasing the set $C$ of signals for player 1 has no influence on the existence of the uniform value, so in the sequel  we will assume w.l.o.g. that: {\it The mapping $(\, c\longrightarrow (\hat{k}(c), \hat{d}(c))\,)$  is a surjection from $C$ to $K \times D$. }
 
 We put $X =\Delta(K)$. An element $u$ in $\Delta_f(X)$ is written $u=\sum_{p \in X} u(p) \delta_p$.  As in the previous section, we use the Wasserstein distance, and the (reverse of ) the Choquet order  on $\Delta(X)$.  $\forall u \in \Delta(X), \forall v \in \Delta(X),  \;\; d(u,v)= \sup_{f:E\to \R, 1-Lip} |u(f)-v(f)|$. And we write  $u\succeq v$   iif for every continuous {\it concave}  real valued mapping $f$ defined on $X$ , $u(f)\geq v(f).$
  
If ${\cal S}$ is a finite set,  we use the norm $\|. \|_1$ on $\R^{\cal S}$.  The set of probability distributions $\Delta({\cal S})$ is viewed\footnote{Notice that if we put $d(s,s')=2$ for any distinct elements of ${\cal S}$, then   for every $p$ and $q$ in $\Delta({\cal S})$ we have  $\sup_{f:{\cal S}\to \R, 1-Lip} |\sum_s p^sf(s)-\sum_s q^sf(s)|= \|p-q\|_1$.}      as a subset of $\R^{\cal S}$.

 \subsubsection{Value of finite games.}\label{subsub6}
 
 As in definition \ref{def9}, we need to consider a large family of finite games. 
 
 \begin{defi} \label{def19} Let        $\theta=\sum_{t\geq 1} \theta_t \delta_t$ be in $\Delta_f(\N^*)$, i.e. $\theta$ is a probability with finite support over positive integers.  For $\pi$ in $\Delta(K \times C \times D)$, the game  $\Gamma_{[\theta]} (\pi)$ is the game with normal form $(\Sigma, {\cal T}, \gamma_{[\theta]}^{\pi})$, where: $$\gamma_{[\theta]} ^\pi (\sigma, \tau)=\E_{\P_{\pi, \sigma, \tau}} \left(  \sum_{t\geq 1}  \theta_t\; g (k_t,i_t,j_t)
\right).$$

\noindent Particular cases:  if $\theta=1/n\; \sum_{t=1}^n \delta_t$, $\Gamma_{[\theta]} (\pi)$ is nothing but $\Gamma_n(\pi)$. 
 
For  $m\geq 0$ and $n\geq 1$,  we denote by $\Gamma_{m,n}(\pi)$ the game  $\Gamma_{[\theta]} (\pi)$ where $\theta=1/n\; \sum_{t=m+1}^{m+n} \delta_t$. The payoff function is written in this case: $\gamma_{m,n} ^\pi (\sigma, \tau)$.  The value of $\Gamma_{m,n}(\pi)$ will  be  denoted  by $v_{m,n}(\pi)$. \end{defi}
Notice that $v_{0,n}$ is just  $v_n$, the value of the $n$-stage game.  The following lemma is true without the hypotheses $HA'$ and $HB'$.

\begin{lem} \label{lem9} For every $\theta\in \Delta_f(\N^*)$ and $\pi \in \Delta(K \times C \times D)$, the game $\Gamma_{[\theta]}(\pi)$ has a value, denoted by $v_{[\theta]}(\pi)$, and both players have optimal strategies. Moreover, $v_{[\theta]}$ is a non expansive mapping from $\Delta(K \times C \times D)$ to $\R$. \end{lem}

\noindent{\bf Proof:} The existence of the value and optimal strategies is standard. Notice that for every $\theta$, $\pi$, and strategy pair $(\sigma, \tau)$:
$$\gamma_{[\theta]} ^{\pi} (\sigma, \tau)=\sum_{k_1,c_1,d_1} \pi (k_1,c_1,d_1) \,\gamma_{[\theta]} ^{\delta_{(k_1,c_1,d_1)}}(\sigma, \tau).$$
 Since  we use $\|.\|_1$,  $v _{[\theta]}$ is 1-Lipschitz.\hfill $\Box$

\begin{defi} \label{def20} We define a mapping $\Psi$ from $\Delta(K \times D)$ to $\Delta_f(X)$ by: for each probability $\pi$ on $K\times D$, $\Psi(\pi)=\sum_{d\in D} \pi(d) \delta_{\pi^d}$, where for each $d$, $\pi^d$ is the conditional probability on $K$ given $d$ issued from $\pi$.\end{defi}   
 
 \begin{nota} \label{not5}  Let  $\pi$ be in  $\Delta(K \times C \times D)$.   We denote by $\bar{\pi}$ the marginal of $\pi$ on $K \times D$, and  denote by $\hat{\pi}$ the probability induced by $\pi$ (or $\bar{\pi}$) on $X$, i.e. $\hat{\pi}=\Psi(\bar{\pi})=\sum_{d \in D} \pi(d) \delta_{\pi^d}$ $\in$ $\Delta_f(X)$, where for each $d$, $\pi^d$ is the conditional probability on $K$ given $d$ issued from $\pi$ (or $\bar{\pi}$).
  
  We also put   $\Delta^E=\{\pi \in \Delta(K \times C \times D), \pi(E)=1\}$, where $E=\{(k,c,d)\in K \times C \times D, \hat{k}(c)=k\; \mbox{ and } \hat{d}(c)=d\}$.\end{nota}

\begin{lem} \label{lem10} Let $\pi$ and $\pi'$ be in $\Delta^E$ such that $\hat{\pi}=\hat{\pi'}$. Then $v_{[\theta]}(\pi)=v_{[\theta]}(\pi')$ for each $\theta$. \end{lem}

\noindent{\bf Proof:} Fix $\theta$ in $\Delta_f(\N^*)$, $\pi$ in $\Delta^E$, and a strategy pair $(\sigma, \tau)$. $$\gamma_{[\theta]} ^{\pi} (\sigma, \tau)=\sum_{d_1} \pi (d_1)\sum_{k_1}\pi(k_1|d_1)  \sum_{c_1} \pi(c_1|d_1,k_1) \,\gamma_{[\theta]} ^{\delta_{(k_1,c_1,d_1)}}(\sigma, \tau).$$
\noindent Since $\pi(E)=1$, player 1 can deduce $d_1$ and $k_1$ from $c_1$, so:
  $$\sup_{\sigma \in \Sigma} \gamma_{[\theta]} ^{\pi} (\sigma, \tau)=\sum_{d_1} \pi (d_1)\sum_{k_1}\pi(k_1|d_1)  \sum_{c_1} \pi(c_1|d_1,k_1) \,\sup_{\sigma \in \Sigma} \gamma_{[\theta]} ^{\delta_{(k_1,c_1,d_1)}}(\sigma, \tau).$$
  \noindent But $\sup_{\sigma \in \Sigma} \gamma_{[\theta]} ^{\delta_{(k_1,c_1,d_1)}}(\sigma, \tau)$ does not depend on $c_1$, so for any fixed $c^*$ in $C$,  $$\sup_{\sigma \in \Sigma} \gamma_{[\theta]} ^{\pi} (\sigma, \tau)=\sum_{d_1} \pi (d_1)\sum_{k_1}\pi(k_1|d_1)    \,\sup_{\sigma \in \Sigma} \gamma_{[\theta]} ^{\delta_{(k_1,c^*,d_1)}}(\sigma, \tau).$$
  Consequently, $\sup_{\sigma \in \Sigma} \gamma_{[\theta]} ^{\pi} (\sigma, \tau)$ only depends on $\hat{\pi}$, $\tau$ and $\theta$, and $v_{[\theta]}(\pi)$ $=\inf_{\tau \in {\cal T}} $ $\sup_{\sigma \in \Sigma} \gamma_{[\theta]} ^{\pi} (\sigma, \tau)$ only depends  on $\hat{\pi}$ and $\theta$. \hfill $\Box$

\vspace{0,5cm}

Remember that we assumed w.l.o.g. that the function  $(\hat{k},\hat{d})$ appearing in hypothesis HA' is surjective.   It will be convenient in the sequel  to use the following notation. 

\begin{nota} \label{not126}  For any $(k,d)$ in $K \times D$, we fix an element $c(k,d)$ in $C$ such that $\hat{k}(c(k,d))=k$ and $\hat{d}(c(k,d))=d$.\end{nota}

$v_{[\theta]}$ has been defined as a mapping from $\Delta(K\times C \times D)$ to $\R$. We  now define   value functions with domain $X$ and $\Delta(X)$. Let $p$ be in $X$. We define $\pi$ in $\Delta^E$ as follows: fix $d^*$ in $D$, and $\pi$ chooses, for each $k$ in $K$, the element $(k,c(k,d^*),d^*)$  with probability $p^k$. Then $\hat{\pi}=\delta_p$, so by the previous lemma $v_{[\theta]}(\pi)$ only depends on $p$. We thus define : 
 $$v_{[\theta]}(p)=v_{[\theta]}(\pi).$$
\noindent With a slight abuse of notation, $v_{[\theta]}$ now also denotes a mapping from $\Delta(K)$ to $\R$. And we have for each $\pi$ in $\Delta^E$:  \begin{eqnarray*}
v_{[\theta]}(\pi)& =&\inf_{\tau \in {\cal T}} \sup_{\sigma \in \Sigma} \gamma_{[\theta]} ^{\pi} (\sigma, \tau),\\
\; &=& \sum_{d_1\in D} \pi (d_1)\; \inf_{\tau \in {\cal T}} \left( \sum_{k_1}\pi(k_1|d_1)    \,\sup_{\sigma \in \Sigma} \gamma_{[\theta]} ^{\delta_{(k_1,c^*,d_1)}}(\sigma, \tau)\right),\\
\; &=&\sum_{d_1\in D} \pi(d_1) v_{[\theta]}(\pi^{d_1}).
\end{eqnarray*}
So we have obtained the following.

 \begin{lem} \label{lem11} $$\forall \pi \in \Delta^E, \forall \theta \in \Delta_f(\N^*),\;\;  v_{[\theta]}(\pi)=\sum_{d_1\in D} \pi(d_1) v_{[\theta]}(\pi^{d_1}).$$\end{lem}
 
\begin{nota} \label{not6}$\;$ $\tilde{v}_{[\theta]}$ denotes the affine extension of $v_{[\theta]}$ on $\Delta(X)$, i.e.: $\forall u \in \Delta(X),  \tilde{v}_{[\theta]}(u)= \int_{p \in X}v_{[\theta]}(p)du(p).$ \end{nota}

From the previous computations, $\tilde{v}_{[\theta]}$ is clearly   linked to  the original value function $v_{[\theta]}$. 
\begin{cla}\label{cla7} $$\forall \pi \in \Delta^E, v_{[\theta]}(\pi)= \tilde{v}_{[\theta]}(\hat{\pi}).$$\end{cla}
  
 So from the knowledge of $v_{[\theta]}$ on $X$, one can deduce its extension  $\tilde{v}_{[\theta]}$ on $\Delta(X)$ and then the original value function $v_{[\theta]}$ on $\Delta^E$. We have, for each $p$ in $\Delta(K)$ (and for any $d^*$ in $D$):
  
   $v_{[\theta]}(p)= \inf_{\tau \in {\cal T}} \left( \sum_{k_1}p^{k_1}\,\sup_{\sigma \in \Sigma} \, \gamma_{[\theta]} ^{\delta_{(k_1,c(k,d^*),d^*)}}(\sigma, \tau)\right).$ So $v_{[\theta]}$ is a concave and non expansive mapping from $\Delta(K)$ to $\R$. Consequently, $\tilde{v}_{[\theta]}$ is a non decreasing and non expansive mapping from $\Delta(X)$ to $\R$.\\

     We finally define : 
     \begin{defi} \label{def20,5} For $\pi$ in $\Delta(K \times C \times D)$, we put:$$\dis v^*(\pi)=\inf_{n\geq 1} \,  \sup_{m\geq 0} \,  v_{m,n}(\pi).$$ \end{defi}

\subsubsection{An auxiliary stochastic game.}\label{subsub7}
   
   We now introduce a stochastic game with complete information, to be played in pure strategies, as in section \ref{secstoch}.
   
   \begin{defi} \label{def21}  Recall that $X=\Delta(K)$.   We put   $A=\Delta(I)^K$ and $B=\Delta(J)$, and    define for every $p$ in $X$, $a$ in $A$ and $b$ in $B$:   \begin{eqnarray*}
    g(p,a,b)& =& \sum_{k\in K} p^k \sum_{i \in I} \sum_{j \in J} a^k(i)b(j)g(k,i,j),\\
    g(p,a)&=&\inf_{b \in B} g(p,a,b),\\
    Q(p,a,b)& =&\sum_{(k,i,j) \in K\times I \times J} p^ka^k(i)b(j) q(k,i,j)\in \Delta(K \times C \times D), \\
   \bar{Q}(p,a)&=&\sum_{(k,i) \in K\times I} p^ka^k(i)\bar{q}(k,i) \in \Delta(K \times D),\\
    l(p,a)   &= &\Psi(\bar{Q}(p,a)).
    \end{eqnarray*}
 $g$ is a mapping  from $X\times A \times B$ to $[0,1]$, 
whereas  $l$ is a mapping from $X \times A$ to $\Delta_f(X)$. 

For $u$ in $\Delta_f(X)$, we write $\hat{\Gamma}(u)$ for the stochastic game $(X,A,B,g,l,u)$ with initial distribution $u$. \end{defi}

   By hypothesis HB', it is easy to see that the marginal of $Q(p,a,b)$ on $\Delta(K\times D)$ does not depend on $b$, and  precisely is $\bar{Q}(p,a)$. $ l(p,a)$ is nothing but  
   
   \noindent   $\sum_{d\in D}\bar{Q}(p,a)(d)\,  \delta_{\bar{Q}(p,a)^d}$, where for each $d$, $\bar{Q}(p,a)^d$ is the conditional probability on $K$ given $d$ issued from $\bar{Q}(p,a)$. Notice also that for every $(p,a,b)$, $Q(p,a,b)$ belongs to the convex set $\Delta^E$. 
   
   Suppose that in the original game $\Gamma(\pi)$, the current state is selected according to $p$, player 1 knows $k$ and plays the mixed action $a^k\in \Delta(I)$, whereas player 2 just knows $p$ and plays the mixed action $b\in \Delta(J)$. Then $g(p,a,b)$ is the (ex-ante) expected payoff for player 1, and $l(p,a)$ is the (ex-ante) distribution of player 2's future belief on the next state. 
 
  We will eventually apply theorem \ref{thm3} to $\hat{\Gamma}(\hat{\pi})$, so we have to check the hypotheses H1 to H7 of section \ref{secstoch}. H1, H2, H3 and H4 are clearly true. We now need the following properties of the mapping $\Psi$.
  
  \begin{lem} \label{lem12} $\Psi $ is concave:$$\forall \pi', \pi'' \in \Delta(K \times D), \forall \lambda \in [0,1],\;\Psi(\lambda \pi' +(1-\lambda) \pi'') \succeq \lambda\Psi( \pi' )+(1-\lambda) \Psi(\pi'').$$
  \end{lem}
  
  \noindent{\bf Proof:} Write $\pi=\lambda \pi' +(1-\lambda) \pi''$. Notice that for each $d$ in $D$, $\pi^d$ $=$ $\frac{1}{\pi(d)} (\lambda \pi'(d) \pi'^d +(1-\lambda) \pi''(d) \pi''^d)$. Let $f$ be a concave continuous mapping from $X$ to $\R$, we have to show that $\Psi(\pi)(f)\geq \lambda  \Psi(\pi')(f)+ (1-\lambda) \Psi(\pi'')(f).$ \begin{eqnarray*}
  \lambda  \Psi(\pi')(f)+ (1-\lambda) \Psi(\pi'')(f) &=& \lambda \sum_d \pi'(d) f(\pi'^d) +(1-\lambda) \sum_{d} \pi''(d) f(\pi''^d),\\
  & = & \sum_d \pi(d) \left( \frac{\lambda \pi'(d)}{\pi(d)} f(\pi'^d) + \frac{(1-\lambda)\pi''(d)}{\pi(d)} f(\pi''(d)\right),\\
  & \leq& \sum_{d} \pi(d) f\left( \frac{\lambda \pi'(d)}{\pi(d)} \pi'^d + \frac{(1-\lambda)\pi''(d)}{\pi(d)} \pi''(d)\right),\\
 & = & \sum_{d} \pi(d) f(\pi^d)=\Psi(\pi)(f).\hspace{4cm}\Box
  \end{eqnarray*}
For any $p$ in $X$, the marginal $\bar{Q}(p,a)$ is affine in $a$, so we obtain that $l(p,a)=\Psi( \bar{Q}(p,a))$ is concave in $a$. Hypothesis H6 will then immediately follow from the next lemma.

  \begin{lem} \label{lem13} $\Psi$ is continuous. \end{lem}
  
  \noindent{\bf Proof:} Let $(\pi_n)_n$ be a   sequence in $\Delta(K \times D)$ converging   for the norm $\|.\|_1$ to $\pi$. It is easy to see that for every $f$ continuous, $\Psi(\pi_n)(f)=\sum_d \pi_n(d) f(\pi_n^d)$ converges as $n$ goes to infinity to $\sum_d \pi(d) f(\pi^d)=\Psi(\pi)(f)$. \hfill $\Box$
  
 \begin{rem}  One can show that $\Psi$ is  Lipschitz, but it is not  1-Lipschitz when $\|.\|_1$ and the Wasserstein distance are used. \end{rem}

 \begin{lem} \label{lem14}  ``Splitting hypothesis" H7. Consider a   convex combination $p=\sum_{s=1}^S \lambda_s p_s$ in $X$,  and   a family of actions $(a_s)_{s \in S}$ in $A^S$. Then there exists $a$ in $A$ such that:
     $$\dis  l(p,a) \succeq \sum_{s\in S} \lambda_s l(p_s,a_s) \; \;\; \rm { and }\it  \;\; \;g(p,a)\geq \sum_{s\in S} \lambda_s   g(p_s,a_s).$$ \end{lem}
     
     \noindent{\bf Proof:} Define $a: K \longrightarrow \Delta(I)$ with the well known splitting procedure of Aumann and Maschler: observe $k$ in $K$ which has been selected  according to $p$, then choose $s$ with probability $\lambda_s p_s^k/p^k$, and finally play $a_s^k$. Formally, put $a^k=\sum_{s\in S} \frac{\lambda_s p_s^k}{p^k}a_s^k \in \Delta(I)$ if $p^k>0$, and define arbitrarily $a^k$ if $p^k=0$. We have:
     \begin{eqnarray*}
     \bar{Q}(p,a) &=& \sum_{k \in K} \sum_{i \in I} p^k a^k(i) \bar{q}(k,i),\\
     &=& \sum_{s\in S} \sum_{k \in K} \sum_{i \in I} \lambda_s p_s^k a_s^k(i) \bar{q}(k,i),\\
     &=&\sum_{s \in S}  \lambda_s \bar{Q}(p_s,a_s)
     \end{eqnarray*}
     \noindent So by concavity of $\Psi$, we have $l(p,a)\succeq \sum_s \lambda_s l(p_s,a_s).$

     Regarding payoffs, we have for each $b$ in $B$, $g(p,a,b)=\sum_s \lambda_s g(p_s,a_s,b)$, so $\inf_{b \in B} g(p,a,b) \geq \sum_{s\in S} \lambda_s \inf_{b \in B} g(p_s,a_s,b).$ \hfill $\Box$
     
Up to now, only H5 remains to be proved.  
     
\subsubsection{The recursive formula.}  \label{subsub8}
We now prove a standard recursive formula for the value functions. As in definition \ref{def9},  for $\theta=\sum_{t \geq 1} \theta_t \delta_t$ we  define $\theta^+$ as the law of $t^*-1$ given that $t^*\geq 2$, so that  $\theta^+=\sum_{t \geq 1} \frac{\theta_{t+1}}{1-\theta_1} \delta_t$  if $\theta_1\neq 1$, and $\theta^+$ is defined arbitrarily    if $\theta_1= 1$.
   
   \begin{pro} \label{pro8}     
 For  $\theta$ in $\Delta_f(\N^*)$ and $p$ in $X$,  \begin{eqnarray*} 
  v_{[\theta]}(p) & = &\max_{a \in A} \min_{b \in B} \left( \theta_1 g(p,a,b)+ (1-\theta_1) \tilde{v}_{[\theta^+]} (l(p,a)) \right),\\
\; & = & \min_{b \in B} \max_{a \in A} \left( \theta_1 g(p,a,b)+ (1-\theta_1) \tilde{v}_{[\theta^+]} (l(p,a)) \right).
 \end{eqnarray*}
For every $\pi$ in $\Delta^E$,  in the game $\Gamma_{[\theta]}(\pi)$ both players have optimal strategies only depending on $\bar{\pi}\in \Delta(K\times D)$. \end{pro}

\noindent{\bf Proof:} By the proof of lemma \ref{lem10}, we know that for any $\tau$ in ${\cal T}$ and fixed $c^*$ in $C$,  $$\sup_{\sigma \in \Sigma} \gamma_{[\theta]} ^{\pi} (\sigma, \tau)=\sum_{d_1} \pi (d_1)\sum_{k_1}\pi(k_1|d_1)    \,\sup_{\sigma \in \Sigma} \gamma_{[\theta]} ^{\delta_{(k_1,c^*,d_1)}}(\sigma, \tau).$$
\noindent And $\tau$ is optimal in $\Gamma_{[\theta]}(\pi)$  if and only if  $\sum_{d_1} \pi (d_1)\sum_{k_1}\pi(k_1|d_1)    \,\sup_{\sigma \in \Sigma} \gamma_{[\theta]} ^{\delta_{(k_1,c^*,d_1)}}(\sigma, \tau)$ $=$ $\sum_{d_1} \pi (d_1) \inf_{\tau' \in {\cal T}}\sum_{k_1}\pi(k_1|d_1)    \,\sup_{\sigma \in \Sigma} \gamma_{[\theta]} ^{\delta_{(k_1,c^*,d_1)}}(\sigma, \tau').$ Hence player 2 has an optimal strategy in $\Gamma_{[\theta]}(\pi)$ that only depends on $\bar{\pi}$.

We now show that player 1 has an optimal strategy in $\Gamma_{[\theta]}(\pi)$ that only depends on $\bar{\pi}$. For every $d_1$ in $D$, fix an optimal strategy $\sigma(d_1)$ of player 1 in the game $\Gamma_{[\theta]}\left(\sum_{k_1\in K} \pi(k_1|d_1) \delta_{(k_1,c(k_1,d_1),d_1)} \right)$. Define now $\sigma$ that plays after each  initial signal $c_1$ exactly what $\sigma(\hat{d}(c_1))$ plays after the initial signal $c(\hat{k}(c_1),\hat{d}(c_1))$. One can check that $\sigma$   is optimal in any game $\Gamma_{[\theta]}(\pi')$, if $\pi'\in \Delta^E$ and $\bar{{\pi'}}=\bar{\pi}$.

We now prove the recursive formula by induction on the greatest element in the support of $\theta$. $v_1(p)=\max_{a \in A} \min_{b \in B}  g(p,a,b)= \min_{b \in B}  \max_{a \in A} g(p,a,b)$ is easy by Sion's theorem. 

  Fix now $n \geq2$, and assume that the proposition is true for every $\theta$ with support included in $\{1,...,n-1\}$. Fix a probability $\theta=\sum_{t=1}^n \theta_t \delta_t$, and notice that $\theta^+$ has a support included  in $\{1,...,n-1\}$. Fix also $p$ in $X$. The equality   $\max_{a \in A} \min_{b \in B}$ $  ( \theta_1 g(p,a,b)+ (1-\theta_1) \tilde{v}_{[\theta^+]} (l(p,a)) )$ $=$ $ \min_{b \in B} \max_{a \in A} $ $ ( \theta_1 g(p,a,b) + (1-\theta_1) \tilde{v}_{[\theta^+]} (l(p,a))  )$ is standard and similar to the proof of proposition \ref{pro2}, so the proof is omitted here. 

By definition, we have $v_{[\theta]}(p)=v_{[\theta]}(\pi)$, where $\pi=\sum_{k} p^k \delta_{(k,c(k,d^*),d^*)} \in \Delta(K\times C \times D)$, and $d^*$ is an arbitrary element of $D$. We thus consider the game $\Gamma_{[\theta]}(\pi)$.  Let  $(\sigma, \tau)$ be  a strategy pair.  Write $a=(a^k)_k=(\sigma_1(c(k,d^*)))_k$ in $A=\Delta(I)^K$ for the action played by player 1 at stage 1, and write $b=\tau_1(d^*)\in B=\Delta(J)$ for the first stage action of player 2. We have:$$\gamma_{[\theta]}^{\pi}(\sigma, \tau)= \theta_1 g(p,a,b) + (1-\theta_1) \sum_{k,i_1,j_1} p^ka^k(i_1)b(j_1) \gamma_{[\theta^+]}^{q(k,i_1,j_1)} (\sigma^+_{c(k,d^*),i_1}, \tau^
+_{d^*,j_1}),$$
\noindent where $\sigma^+_{c(k,d^*),i_1}$ and  $\tau^
+_{d^*,j_1}$ are continuation strategies.

$a$ being fixed, we can choose $\sigma^+$ an optimal strategy for player 1 in the game $\Gamma_{[\theta^+]}(Q(p,a,j_1))$, and this choice can be made independently of $j_1$ since the marginal $\bar{Q}(p,a,j_1)$ does not depend on $j_1$. We have for every $j_1$ in $J$: $ \sum_{k,i_1} p^ka^k(i_1) \gamma_{[\theta^+]}^{q(k,i_1,j_1)} (\sigma^+, \tau^+_{d^*,j_1})$ $=\gamma_{[\theta^+]}^{Q(p,a,j_1)} (\sigma^+, \tau^+_{d^*,j_1})$ $\geq v_{[\theta^+]}(Q(p,a,j_1))$ $=$ 

\noindent $\tilde{v}_{[\theta^+]}(l(p,a)).$ By playing $a$ at stage 1, and then according to $\sigma^+$ (for any first signal  $c(k,d^*)$ and   first action $i_1$),  player 1 can thus guarantee: $\inf_b \theta_1 g(p,a,b) + (1-\theta_1) \tilde{v}_{[\theta^+]}(l(p,a)).$ So $v_{[\theta]}(p)\geq \max_{a \in A} \min_{b \in B} ( \theta_1 g(p,a,b)+$ $ (1-\theta_1) \tilde{v}_{[\theta^+]} (l(p,a)) )$.

We finally show that player 2 can defend $\max_{a \in A} \min_{b \in B} ( \theta_1 g(p,a,b)$ $+ (1-\theta_1) \tilde{v}_{[\theta^+]} (l(p,a)) )$ in $\Gamma_{[\theta]}(\pi)$. Fix  a strategy $\sigma$ of player 1. $a=(\sigma_1(c(k,d^*)))_k$   being fixed, choose $b$ in $B$ achieving $\inf_b g(p,a,b)$. We also   choose $\tau^+$ in ${\cal T}$ an optimal strategy for player 2 in the game $\Gamma_{[\theta^+]}(Q(p,a,j_1))$, and this choice can be made independently of $j_1$. $\tau^+$ now  being   fixed, notice that there exists a   strategy $\sigma'$ of player 1 which is a best reply to $\tau^+$ in any game $\Gamma_{[\theta^+]}(\pi')$, with $\pi'\in \Delta^E$. We now have:
\begin{eqnarray*}
\gamma_{[\theta]}^{\pi}(\sigma, \tau) &= &\theta_1 g(p,a,b) + (1-\theta_1) \sum_{j_1} b(j_1) \sum_{k,i_1} p^ka^k(i_1)\gamma_{[\theta^+]}^{q(k,i_1,j_1)} (\sigma^+_{c(k,d^*),i_1}, \tau^
+),\\
 & \leq & \theta_1 g(p,a) + (1-\theta_1) \sum_{j_1} b(j_1) \sum_{k,i_1} p^k a^k(i_1)\gamma_{[\theta^+]}^{q(k,i_1,j_1)} (\sigma', \tau^+),\\
 & = & \theta_1 g(p,a) + (1-\theta_1) \sum_{j_1} b(j_1) \gamma_{[\theta^+]}^{Q(p,a,j_1)} (\sigma', \tau^
+),\\
 & = & \theta_1 g(p,a) + (1-\theta_1) \sum_{j_1} b(j_1) v_{[\theta^+]}(l(p,a)),\\
 & \leq &   \sup_{a \in A}  \left( \theta_1 g(p,a)+ (1-\theta_1) \tilde{v}_{[\theta^+]} (l(p,a)) \right)
\end{eqnarray*}
So  $v_{[\theta]}(p)\leq \max_{a \in A} \min_{b \in B} \left( \theta_1 g(p,a,b)+ (1-\theta_1) \tilde{v}_{[\theta^+]} (l(p,a)) \right)$.  \hfill $\Box$

  \subsubsection{Player 2 can guarantee $v^*(\pi)$ in $\Gamma(\pi)$.} \label{subsub9}
  
  We first have an analog of lemma \ref{lem5}. Recall that a strategy for player 2 is a sequence  $\tau=(\tau_1,\tau_2,...,\tau_t,...)$, where for each $t$, $\tau_t$ is a mapping from $(D \times J)^{t-1}\times D$ to $\Delta(J)$. 
  
  \begin{lem} \label{lem15} For every $\pi\in \Delta^E$, $n\geq 1$ and $m\geq 0$, $\forall \tau_1,...,\tau_m,$ $\exists \tau_{m+1},...,\tau_{m+n}$ such that any strategy of player 2 starting by $\tau_1,...,\tau_{m},...,\tau_{m+n}$ is optimal in the game $\Gamma_{m,n}(\pi)$.\end{lem}
  
  \noindent{\bf Proof:} Fix $\pi$, $m$, $n$, and $\tau_1$, ..., $\tau_m$, with $\tau_t:(D \times J)^{t-1}\times D\longrightarrow    \Delta(J)$ for every $t\leq m$. Define ${\cal T}^{\dag}$ as the set of strategies of player 2 that start with $\tau_1$,...,$\tau_m$. Let us now consider the zero-sum game $\Gamma^{\dag}_{m,n}(\pi)$ with strategy set $\Sigma $ for player 1, ${\cal T}^{\dag}$ for player 2, and payoff function (the restriction of) $v_{m,n}$. Stages greater than $m+n$ do not care, so $\Gamma^{\dag}_{m,n}(\pi)$  can be seen as the mixed extension of a finite game, where nature plays $\tau_1$, ..., $\tau_m$ instead of player 2 for the first $m$ stages. Consequently,    $\Gamma^{\dag}_{m,n}(\pi)$ has a value which we denote by $v^{\dag}_{m,n}(\pi)$. Clearly, $v^{\dag}_{m,n}(\pi)\geq  v_{m,n}(\pi)$. Now, for any strategy $\sigma$ of player 1, it is easy to construct, by the recursive formula of proposition \ref{pro8}, a strategy $\tau$ that defends $v_{m,n}(\pi)$ against $\sigma$. So  $v^{\dag}_{m,n}(\pi)=v_{m,n}(\pi)$, and considering an optimal strategy of player 2 in $\Gamma^{\dag}_{m,n}(\pi)$ concludes the proof. \hfill $\Box$
  
  \begin{pro} \label{pro10} For each $\pi$ in $\Delta^E$, player 2 can guarantee $v^*(\pi)$ in the game $\Gamma(\pi)$.\end{pro}
  
  \noindent {\bf Proof:} Divide the set of stages into consecutive blocks $B^1$, $B^2$,..., $B^m$ of equal length $n$. Define  a strategy $\tau$ of player 2 as follows. At block $B^1$, pick $\tau_1, \tau_2,..., \tau_n$ in order to get an optimal strategy in $\Gamma_{0,n}(\pi)$. At block $B^2$,    use lemma \ref{lem15} to construct   $\tau_{n+1},...,\tau_{2n}$ and  get an optimal strategy also in $\Gamma_{n,n}(\pi)$, etc... At block $B^{m+1}$, given $\tau_1$,...,$\tau_{nm}$, use lemma \ref{lem15} to define $\tau_{nm+1}$,..., $\tau_{n(m+1)}$ to get an optimal strategy in $\Gamma_{nm,n}(\pi)$. 
  
  For any $\sigma$ and $M\geq 1$, we have:
   
  $\E_{\P_{\pi, \sigma, \tau}} \left( \frac{1}{Mn} \sum_{t=1}^{Mn}g(k_t,i_t,j_t)\right) \leq \frac{1}{M} \sum_{m=0}^{M-1} \tilde{v}_{mn,n}(\pi) \leq \sup_m  \tilde{v}_{m,n}(\pi).$   So player 2 can guarantee  $v^*(\pi)$. \hfill $\Box$
  
  \vspace{0,5cm}
  
  We  now have the following inequality chain, for $\pi$ in $\Delta^E$:
 $$\boxed{\underline{v}(\pi) \leq v^-(\pi)\leq  v^+(\pi) \leq \overline{v}(\pi)\leq v^*(\pi).}$$

\subsubsection{Player 1 can guarantee $v^*(\pi)$ via the auxiliary game.} \label{subsub10} 

For  $f$ continuous and $\alpha$ in $[0,1]$, we defined in notation \ref{nota2} the mapping  $\Phi(\alpha,f)$ as follows: $\forall p \in X, $

$\;\; \Phi(\alpha,f)(p)= \sup_{a \in A} \inf_{b \in B} \;\left( \;  \alpha \; g(p,a,b) +(1-\alpha) \; \tilde{f} (l(p,a))\; \right).$ 

\noindent We now simply define the following subset of mappings from $X$ to $\R$:   $${\cal D}=\{v_{[\theta]}, \; \theta \in \Delta_f(\N^*)\}.$$ 
By the recursive formula, we obtain that ${\cal D}$ is stable under $\Phi$:    $\forall f \in {\cal D}, \forall \alpha \in [0,1]$,   $ \Phi(\alpha,f) \in {\cal D}$. Since $v_1=\Phi(1,0)\in {\cal D}$ and all elements of ${\cal D}$ are non expansive, the  hypothesis H5 of section \ref{secstoch} is satisfied.  

\begin{pro} \label{pro9} $\;$

a)  For every $u$ in $\Delta_f(X)$, the auxiliary game $\hat{\Gamma}(u)=(X,A,B,g,l,u)$ satisfies   the hypotheses H1,..., H7 of theorem \ref{thm3}.

b) For every $\theta$ in $\Delta_f(\N^*)$ and $u$ in $\Delta_f(X)$, the  auxiliary game  $\hat{\Gamma}_{[\theta]}(u)$  defined in definition \ref{def9} has a value which corresponds to $\tilde{v}_{[\theta]}(u)$, as defined in notation \ref{not6}. 

c) For $\pi$ in $\Delta^E$, anything that can be guaranteed by player 1 with Markov strategies in $\hat{\Gamma}_{[\theta]}(\hat{\pi})$  can be guaranteed by player 1 in  the original game $\Gamma(\pi)$. \end{pro}

\noindent{ \bf Proof:} a) H1,...,H7 have been proved. b) The equality between the value functions of the original game and of the auxiliary  game comes from propositions \ref{pro2} and \ref{pro8}. c) A markov strategy of player 1 is a sequence  $(\sigma_t)_{t\geq 1}$, where for each $t$ $\sigma_t$ is a mapping from $X$ to $A$ giving the action to be played on stage $t$ depending on the current state in $X=\Delta(K)$. It induces a probability distribution on $(X \times A)^{\infty}$, regardless of player 2's actions (see subsection \ref{subsub3}). Any such strategy can be mimicked  by player 1 in the original game, since this player can compute at each stage the belief of player 2 on the state in $K$. \hfill $\Box$

\vspace{0,5cm}

Notice that the analog of point c) is not true regarding player 2, because in the original game this player can compute posterior beliefs  on $K$ only if he knows player 1's strategy. 

\vspace{0,5cm}

We can now conclude the proof of theorem \ref{thm5}. Fix $\pi$ in $\Delta^E$, and put $u=\hat{\pi} \in \Delta_f(X)$. For every $\theta$, we have $v_{[\theta]}(\pi)=\tilde{v}_{[\theta]}(u)$ by proposition \ref{pro9} b) and claim \ref{cla7}. So $v^*(\pi)=\inf_n \sup_m v_{m,n}(\pi)=\inf_n \sup_m \tilde{v}_{m,n}(u)=v^*(u)$. By theorem \ref{thm3}, $(\tilde{v}_n(u))_n$ converges to  $v^*(u)$,  and player 1 can guarantee $v^*(u)$ in $\hat{\Gamma}(u)$ with a Markov strategy. So we obtain that  $({v}_n(\pi))_n$ converges to  $v^*(\pi)$, and player 1 can guarantee this quantity in the original game $\Gamma(\pi)$. Finally  $\Gamma(\pi)$  has a uniform value which is:    $$\boxed{\underline{v}(\pi)= v^-(\pi)= v^+(\pi) = \overline{v}(\pi)=v^*(\pi).}$$ 

\centerline{$\Box$ \hspace{0,2cm} $\Box$ \hspace{0,2cm} $\Box$ \hspace{0,2cm} $\Box$}

\subsection{Comments and consequences} \label{subcomcon}

\subsubsection{Byproducts of the proof.}\label{subsub11}

 The proof  of theorem \ref{thm5} shows, under the very same hypotheses HA' and HB', more than the existence of the uniform value.  

 In particular, the application of theorem \ref{thm3} to the auxiliary game gives:  
\begin{eqnarray*}
v^*(\pi)& = &\inf_{n\geq1} \sup_{m\geq 0}  {v}_{m,n}(\pi)   =   \sup_{m\geq 0}  \inf_{n\geq1} {v}_{m,n}(\pi),\\
 \;& = &\inf_{n\geq1} \sup_{m\geq 0} w_{m,n}(\pi)= 
\sup_{m\geq 0}  \inf_{n\geq1} w_{m,n}(\pi).
 \end{eqnarray*}
\noindent where $  {v}_{m,n}(\pi)$ is defined in definition \ref{def19}, and $w_{m,n}(\pi)=\inf_{\theta \in \Delta(\{1,...,n\})}   v_{[\theta^{m,n}]}(\pi)$ (see definition \ref{def14} and corollary \ref{cor2}).
 
And   $(v_n)_n$ uniformly converges to $v^*$ on $\Delta^E$. 

Concerning $\varepsilon$-optimal strategies, we have seen that player 1 can guarantee $v^*(u)$ with Markov strategies, i.e. with strategies that play at each stage a mixed action  determined by player 2's current belief on the current state of nature in $K$. 

Regarding player 2, one can strengthen the construction of proposition \ref{pro10} and show  as in remark \ref{rem5}   that player 2 has 0-optimal strategies in the game $\Gamma(\pi)$. Notice that our proof does not tell if player 1 has 0-optimal strategies in the game $\Gamma(\pi)$ (see example 5.7. in Renault, 2007).
 \hfill $\Box$

\vspace{0,5cm}

\subsubsection{HA'  or HB' can not be withdrawn in theorem \ref{thm5}.}

An example of a game satisfying HA' and having no uniform value is given in  Sorin, 1984.  It is a particular case of a stochastic game with incomplete information, where after each stage the players perfectly observe the actions just chosen. (There are two possible  stochastic games of  ``Big Match" type, and player 1 only knows which one is being played.)

An example of a game satisfying HB' and having no uniform value is given in  Sorin and Zamir, 1985. It belongs to a class of games called  repeated games   with incomplete information on one and a half side: at each stage, both players will play the same matrix game. Player 1 initially receives a signal which tells him which matrix game will be played, but does not know the initial signal of player 2, so can not deduce from his signal the belief of player 2 on the selected matrix game. The transition function $q$ is particularly simple there: $q(k,i,j)$ is the Dirac measure on $(k, j, i)$. 

\vspace{0,5cm}

\subsubsection{Applications.} \label{subsecappli}

Consider   the following model of repeated games with standard monitoring. There are:  a finite set of states $K$, an initial probability $p$ on $\Delta(K)$, finite action sets $I$ and $J$, for each state $k$ a   payoff matrix $(G^k(i,j))_{(i,j)}$,  and for each state $k$ in $K$ and action $i$ in $I$     there is a probability  $l(k,i)\in \Delta(K)$. At stage 1, a state $k_1$ is selected according to $p$ and told to player 1 only. Then simultaneously player 1 chooses $i_1$ in $I$ and player 2 chooses $j_1$ in $J$. The payoff for player 1 is $G^{k_1}(i_1,j_1)$, and $(i_1,j_1)$ is publicly announced. At stage $t\geq 2$, $k_t$ is selected according to $l(k_{t-1},i_{t-1})$ and told to player 1 only. Then the players choose $i_t$ and $j_t$. The stage payoff for player 1 is $G^{k_t}(i_t,j_t)$,  $(i_t,j_t)$ is publicly announced, and the play goes to stage $t+1$. 

This model is a generalization of the model of Markov chain repeated games with lack of information on one side introduced in Renault, 2006.  Here, player 1 is not only informed of the sequence of actions, but also he can influence the state process. Studying this model  has lead to the present paper, and some ideas developed here already come from Renault 2006. It also contains stochastic games  with a single controller and incomplete information on the side of his opponent, as studied in Rosenberg {\it et al.}, 2004. So the present paper  generalizes both theorem 2.3 in Renault 2006, and theorem 6 in Rosenberg {\it et al.} 2004, and as a consequence it also generalizes the original existence result of Aumann and Maschler (1995) for the value of  (non stochastic)  repeated games with incomplete information on one side and perfect monitoring. 

Notice that our result does not imply the existence of the value for models when player 1 receives signals without having a perfect knowledge of the belief of player 2 on the state (see Aumann Maschler 1995, or Zamir, 1992 for repeated with lack of information on one side, or Neyman 2008 for Markov chain repeated games with lack of information on one side). When the state is uncontrolled,  more flexibility on the signalling structure can be allowed.

\subsubsection{Open problems.}

1. We have seen that   hypotheses HA'  and HB' can not be withdrawn in theorem \ref{thm5}. However, strengthening HB' into HB may allow to weaken   HA' into the following hypothesis.\\

\noindent{ Hypothesis HA''}: Player 1 is more informed than player 2, i.e. there exists a mapping    $\hat{d}:C \longrightarrow D$ such that:  if  $E$ denotes $\{(k,c,d)\in K \times C \times D,  \hat{d}(c)=d\}$, we have: $\pi(E)=1, \; {\rm{and}} \; \; q(k,i,j)(E)=1, \; \forall (k,i,j)\in K \times I \times J.$\\

If we only assume that HA'' and HB' hold, it may be the case that player 2  controls player 1's signal, hence in some sense player 2 may ``manipulate"  player 1's knowledge of the  state.    So our proof   does not apply here, and in our opinion most likely the value may  fail to exist.

The situation is different if we assume  HA'' and HB. Player 1 always have more information than player 2 about  the  state,   but   the set $X=\Delta(K)$ is not sufficient to characterize, after each stage, the difference of information from  player 1 to   player 2. The natural state space here may  rather be the set $\{(u,v)\in \Delta_f(X) \times \Delta_f(X), u \preceq v\}$. Does the uniform value exist  in this case ?\\

2. In general, recall that lemma \ref{lem9} is true without hypotheses, so in particular the $n$ stage values $v_n$ always exist. There is no known example of a zero-sum repeated game (defined with finite data exactly as in subsection \ref{subsub37})  where $\lim_n v_n $ does not exist.\\

3. Assume that player 1 always has more information than player 2, i.e. that player 1 can deduce from his signal both the signal and the action of player 2. This  is the case, e.g., when HA holds. It has been conjectured by Mertens, Sorin and Zamir (see Sorin, 2002, 6.5.8. p.147, or Mertens et al, 1994,  Part C,  p. 451) that  for such repeated games, the limit of  $(v_n(\pi))_n$ exists and can be guaranteed by player 1 in $\Gamma(\pi)$.

The approach used here might help  to prove the conjecture. An important step would be to obtain an analog of the  result on dynamic programming (Renault, 2007) for two-player stochastic games with deterministic transitions and action-independent payoffs. More precisely,  let $Z$ be  a state space,  $A$ and $B$ be   action sets, $r$ be a payoff function from $Z$ to $[0,1]$, and $f$ be a transition from $Z\times A\times B$ to  $Z$. At each stage, if the current state is $z$ simultaneously player 1 chooses $a$ and player 2 chooses $b$.  Player 1's payoff is $r(z)$, and $(a,b)$  and the new state   $f(z,a,b)$ are publicly announced. $Z$ being a precompact metric space, can we find ``nice" uniform equicontinuity conditions on  some auxiliary value functions that will ensure the existence\footnote{Mertens and Neyman (1981)   already provided   conditions on the discounted values which guarantee    the existence of the uniform value in stochastic games, but it seems difficult  to apply them here.} of the uniform value ?  

 \vspace{1cm}
 
  \noindent \large \bf  Acknowledgements. \rm \small 
I thank Guillaume Carlier for stimulating  discussions about Choquet's order and the Wasserstein distance.

The work of Jérôme Renault is currently supported by the GIS X-HEC-ENSAE in Decision Sciences. Most of the present  work was done while the author was at Ceremade, University Paris-Dauphine.
It has been  partly supported by the   French Agence Nationale de la Recherche (ANR), undergrants  ATLAS     and  Croyances, and the ``Chaire de la Fondation du Risque", Dauphine-ENSAE-Groupama : Les particuliers face aux risques.

 \vspace{1cm}

   \vspace{2cm}
   
\noindent \large \bf References. \rm \small

\vspace{0,5cm}

\noindent  Aumann, R.J. and M. Maschler (1995): Repeated games with
 incomplete information.  With the collaboration of R. Stearns.
 {Cambridge, MA: MIT Press.}\\

\noindent  Ash, R.B. (1972): Real Analysis and Probability, Probability and Mathematical Statistics, Academic Press.\\

\noindent Coulomb, J.M. (2003): Games with a recursive structure. based on a lecture of J-F. Mertens. Chapter 28, Stochastic Games and Applications, A. Neyman and S. Sorin eds, Kluwer Academic Publishers.\\

\noindent 
Doob, J.L. (1994):
  Measure Theory.  Springer-Verlag, New-York.\\

\noindent Malliavin, P. (1995)  :   Integration and probability. Springer-Verlag, New-York.  \\

\noindent 
Mertens, J.F. (1986)
\newblock Repeated games. 
\newblock {\em Proceedings of the International Congress of Mathematicians, Berkeley 1986}, 1528--1577. {\it American Mathematical Society}, 1987.\\

\noindent Mertens, J-F. and A. Neyman (1981):  Stochastic games. International Journal of Game Theory, 1, 39-64.\\

\noindent  Mertens, J.F., S. Sorin and S. Zamir (1994): {Repeated Games}.   Parts A, B and C.  CORE Discussion Papers  9420, 9421 and 9422.\\

\noindent Meyer, P.A. (1966):  Probabilit\' es et potentiel. Hermann.\\

\noindent Neyman, A.  (2008):  Existence of Optimal Strategies in Markov Games with Incomplete Information. International Journal of Game Theory, 581-596.\\

\noindent Renault, J. (2006): The value of Markov chain games with lack of information on one side. Mathematics of Operations Research, 31, 490-512.\\

 \noindent Renault, J. (2007): Uniform value in Dynamic Programming. revised version 2009:  http://fr.arxiv.org/abs/0803.2758
 \\

\noindent Rosenberg, D., Solan, E. and N. Vieille (2004): Stochastic games with a single controller and incomplete information. SIAM Journal on Control and Optimization, 43, 86-110. \\

\noindent  Sorin, S.  (1984):  Big match with lack of information on one side (Part I), International Journal of Game Theory, 13, 201-255.  \\

\noindent  Sorin, S.  (2002): A first course on Zero-Sum repeated games. Mathématiques et Applications, SMAI,  Springer.\\

\noindent  Sorin, S. and S. Zamir  (1985): A 2-person game with lack of information on 1 and 1/2 sides. Mathematics of Operations Research, 10, 17-23. \\

\noindent  Zamir, S. (1992): Repeated Games of Incomplete Information: Zero-Sum. {Handbook of
Game Theory with Economic Applications}, edited by Aumann and Hart. Vol I, chapter 6.

 \end{document}